\documentclass{amsart}
\usepackage{amsmath}
\usepackage{amssymb}
\usepackage{amsfonts}
\usepackage{MnSymbol}
\usepackage{graphics} 
\usepackage{epsfig} 
\usepackage{psfrag}
\usepackage{pdfsync}
\usepackage[usenames]{color}
\usepackage{cancel}

\definecolor{arancio}{rgb}{0.90,0.50,0.20}

\definecolor{blu}{rgb}{0.,0.,1.}

\definecolor{siena}{rgb}{0.99,0.39,0.0}

\definecolor{pavone}{rgb}{0.00,0.00,0.63}

\definecolor{malva}{rgb}{0.10,0.50,0.50}

\definecolor{rosso}{rgb}{1.,0.,0.}

\definecolor{geranio}{rgb}{0.90,0.00,0.20}

\definecolor{cerulean}{rgb}
{0.0, 0.48, 0.65}

\newtheorem{theorem}{Theorem}[section]
\newtheorem{corollary}[theorem]{Corollary}

\newtheorem{lem}[theorem]{Lemma}
\newtheorem{prop}[theorem]{Proposition}

\theoremstyle{definition}
\newtheorem{definition}{Definition}[section]
\newtheorem{remark}{Remark}[section]

\newcommand{\ep}{\varepsilon}
\newcommand{\N}{\mathbb{N}}
\newcommand{\R}{\mathbb{R}}

\date{\today}

\newcommand{\bcl}{\begin{center}}
\newcommand{\ecl}{\end{center}}
\newcommand{\brl}{\begin{right}}
\newcommand{\erl}{\end{right}}
\newcommand{\ben}{\begin{enumerate}}
\newcommand{\barr}{\begin{array}}
\newcommand{\earr}{\end{array}}
\newcommand{\btab}{\begin{tabular}}
\newcommand{\etab}{\end{tabular}}
\newcommand{\bdoc}{\begin{document}}
\newcommand{\edoc}{\end{document}}
\newcommand{\beqy}{\begin{eqnarray}}

\newcommand{\beq}{\begin{equation}}
\newcommand{\beqi}{\begin{eqnarray*}}
\newcommand{\bitem}{\begin{itemize}}
\newcommand{\brem}{\begin{remark}}
\newcommand{\erem}{\end{remark}}
\newcommand{\eitem}{\end{itemize}}
\newcommand{\nln}{\newline}
\newcommand{\newt}{\newtheorem}

\renewcommand{\a }{\alpha }
\renewcommand{\b }{\beta }
\newcommand{\g }{\gamma}
\newcommand{\G }{\Gamma }
\renewcommand{\d }{\delta }

\newcommand{\D }{\Delta }
\newcommand{\e }{\epsilon }
\newcommand{\z }{\zeta }
\renewcommand{\l }{\lambda }
\renewcommand{\L }{\Lambda }
\newcommand{\m }{\mu }
\newcommand{\n }{\tau }
\renewcommand{\r }{\rho }
\newcommand{\s }{\sigma }
\newcommand{\Sig }{\Sigma }
\renewcommand{\t }{\tau }
\newcommand{\var }{H }
\renewcommand{\o }{\omega }
\renewcommand{\O }{\Omega }

\newcommand{\supp}{\text{\rm supp}\,}
\newcommand{\sgn}{\text{\rm sgn}\,}

\title[Hamilton-Jacobi equation]
{Discontinuous viscosity solutions \\ of first order Hamilton-Jacobi equations}

\author[Bertsch]{Michiel Bertsch}
\address{Dipartimento di Matematica, Universit\`a di Roma "Tor Vergata", 
Via della Ricerca Scientifica, 00133 Roma, Italy \\ and
Istituto per le Applicazioni del Calcolo "M. Picone", CNR, Roma, Italy} 
\email{bertsch.michiel@gmail.com}
\author[Smarrazzo]{Flavia Smarrazzo}
\address{Universit\`a Campus Bio-Medico di Roma\\ Via Alvaro del Portillo 21, 00128 Roma, Italy}
\thanks{}
\email{flavia.smarrazzo@gmail.com}
\author[Terracina]{Andrea Terracina}
\address{Dipartimento di Matematica "G. Castelnuovo", Universit\`a ''Sapienza'' di Roma\\ P.le A. Moro 5, I-00185 Roma, Italy}
\email{terracina@mat.uniroma1.it}
\author[Tesei]{Alberto Tesei}
\address{Dipartimento di Matematica "G. Castelnuovo", Universit\`a ''Sapienza'' di Roma\\ P.le A. Moro 5, I-00185 Roma, Italy, and 
Istituto per le Applicazioni del Calcolo "M. Picone", CNR, Roma, Italy}
\email{albertotesei@gmail.com}



\date{}

\begin{document}

\bibliographystyle{h-elsevier2}


\begin{abstract}
We consider the simplest example of a time-dependent first order Hamilton-Jacobi equation, i.e.\ in one space dimension and with a bounded and Lipschitz continuous Hamiltonian which  only depends on the spatial derivative. We show that if the initial function has a finite number of jump discontinuities, the corresponding {\em discontinuous} viscosity solution of the corresponding Cauchy problem on the real line is unique. Uniqueness follows from a comparison theorem for semicontinuous viscosity sub- and supersolutions, using the {\it barrier} effect  of spatial discontinuities of a solution. Discontinuous solutions are defined in the spirit of Ishii, but 
semi-continuous envelopes are defined through {\it essential} limits, a definition which is shown to be compatible with Perron's method for existence. We also describe some properties of the evolution of jump discontinuities.
As a by-product we obtain several results on singular Neumann problems.
\end{abstract}

\maketitle

\section{Introduction}\label{intro}
After the introduction of continuous viscosity solutions of first order Hamilton-Jacobi (HJ) equations (\cite{CL1, CL2}), it was readily understood that the basic concepts and comparison results of the theory could be extended to the case of {\em semicontinuous} viscosity  sub- and supersolutions.  A systematic study of discontinuities of both the Hamiltonian itself and solutions of HJ equations, started in \cite{I1,I2}, is an important issue since discontinuous solutions arise in many applications (e.g., optimal control problems,  differential game theory; see \cite{Ba1, CS2} and references therein).

Existence of possibly discontinuous viscosity solutions was proven in \cite{I2} by Perron's method, but 
uniqueness of such solutions remained unclear since in general the classical comparison result for semicontinuous 
viscosity sub- and supersolutions does not imply uniqueness of viscosity solutions - apart from the trivial 
case of continuous  data, when the unique viscosity solution is continuous (\cite{Bl}). In addition, the literature contains interesting examples of nonuniqueness of viscosity solutions for non-convex Hamiltonians with explicit space and/or time dependence and discontinuous initial data (\cite{Ba1, BSS, GS}). 

Motivated by these difficulties, several different notions  of discontinuous solutions of HJ equations have been proposed (\cite{BG, BP, BJ, CS1,GS,Su}) to prove existence, comparison and uniqueness results under various assumptions (for instance, if the Hamiltonian is convex). Although interesting in their own right and for specific applications (e.g.\ to control problems), the relationship between different notions of solution have been elucidated only in some specific cases (see \cite{CS2,GS}).

In the present paper we consider discontinuous viscosity solutions, defined in the spirit of \cite{I2}, in the simplest 
example of a time-dependent HJ equation, i.e.\ the one-dimensional  equation
\begin{equation}\label{hje}
u_t+ H(u_x)=0\,
\end{equation}
with a bounded and uniformly Lipschitz continuous Hamiltonian:
$$
H\in  W^{1,\infty}(\R)\,.
\leqno{(H_1)}
$$
We do not require convexity conditions on $H$, nor the existence of $\lim\limits_{s\to\pm\infty} H(s)$. The boundedness assumption for $H$ is suggested by a model for ion etching (\cite{F, R1,R2}).

Given $u_0\in  L^\infty_{{\rm loc}}(\R)$, we consider in particular the Cauchy problem
$$
\begin{cases}
u_t+ H(u_x)=0&\text{in }\R\times (0,T)\\
u=u_0&\text{in } \R\times\{0\}\,.
\end{cases}
\leqno{(CP)}
$$
As in Ishii's paper (\cite{I2}), the definition of (possibly discontinuous) viscosity sub- and supersolutions is based on 
the concept of semicontinuous envelopes of functions, but we use {\it essential} limits to define semicontinuous envelopes (see Section \ref{prelim}). The use of essential limits excludes some unnatural and artificial examples of admissible viscosity solutions (see for example to the discussion on page 27  of \cite{CIL}). As in \cite{I2}, existence of viscosity solutions of $(CP)$ is based on Perron's method, which can be adapted to our definitions (see Section \ref{perron}). Also the basic comparison result continues to hold for discontinuous viscosity sub- and supersolutions (Section \ref{compre}); since, by $(H_1)$, the speed of propagation is bounded by the Lipschitz constant $\|H'\|_\infty$ (\cite{CL2, I3}), the comparison result is formulated locally.

Since $H$ is bounded, it is rather intuitive that  solutions are Lipschitz continuous with respect to time, and that spatial jump discontinuities of the solution do not disappear instantaneously. In Section \ref{Sec reg} we show that this is indeed the case, and in addition we prove  that spatial jump discontinuities are non-increasing in time and satisfy an explicit decay rate if $\limsup_{s\to \pm\infty}H(s)\!>\! \liminf_{s\to \pm\infty}H(s)$.

The main purpose of the paper is to show that the viscosity solution of $(CP)$ is unique if $u_0$ has a finite number of jump discontinuities: 
$$
\begin{cases}
u_0\in  L^\infty_{{\rm loc}}(\R),\quad
\mbox{$u_0$ is piecewise continuous in $\R$} 
\\
\text{and has a finite number of jump discontinuities.}\end{cases}. \leqno{(H_2)}
$$
If $u_0$ is continuous, the comparison result implies uniqueness of the viscosity solution. In Section ~\ref{unire} we present the proof of  uniqueness if $u_0$ has jump discontinuities.

We briefly describe our approach. By the boundedness of $H$, a jump discontinuity of $u_0$ at a certain point $c$ does not disappear instantaneously. On the other hand it is known (\cite{E}) that spatial discontinuities of a solution produce a {\it barrier effect}\,: if a viscosity solution $u$ has a spatial jump discontinuity at $c\in \R$ for $t\in [0,\tau)$, the evolution of $u$ in $(-\infty,c)\times (0,\tau)$ is independent of that in $(c,\infty)\times (0,\tau)$. More precisely (see Lemma \ref{prerem}), if the jump $u(c^+,t)-u(c^-,t)$  is positive for $t\in(0,\t)$, then $u$ satisfies on either side of $y$ the {\em singular Neumann problems}
$$
\begin{cases}
u_t+ H(u_x)=0&\text{in }(-\infty,c)\!\times\! (0,\tau)\\
u_x(c,t)\!=\!\infty 
&\text{for }0<t<\tau\\
u(x,0)=u_0(x)&\text{for }x<c\,,
\end{cases}
\qquad \ 
\begin{cases}
u_t+ H(u_x)=0&\text{in }(c,\infty)\!\times \!(0,\tau)\\
u_x(c,t)\!=\!\infty 
&\text{for }0<t<\tau\\
u(x,0)=u_0(x)&\text{for }x>c.
\end{cases}
$$
Similarly if $u(c^+,t)-u(c^-,t)<0$, with $u_x(c,t)=\infty$ replaced by $u_x(c,t)=-\infty$.

This leads to the following procedure to prove uniqueness of a discontinuous viscosity solution of $(CP)$. Let $x_1,\ldots,x_p$ be points where $u_0$ has a jump discontinuity, let $u$ and $v$ be two viscosity solutions of $(CP)$  and let $\tau_1\in (0,T]$ be the maximal time for which the jump continuities of $u$ and $v$ at $x_k$ persist in $[0,\tau_1)$ for all $k=1,\ldots,p$.  We define the intervals $I_0=(-\infty,x_1)$, $I_k=(x_k,x_{k-1})$ for $k=1,\ldots, p-1$ and $I_p=(x_p,\infty)$. Then the restrictions of $u$ and $v$ to $I_k\times (0,\tau_1)$ solve the same singular Neumann problem in $I_k\times (0,\tau_1)$ with continuous initial data. Also this singular problem satisfies a comparison principle for viscosity sub- and supersolutions, and since the initial function restricted to $I_k$ is
continuous, this implies that $u=v$ a.e.~in $I_k\times (0,\tau_1)$ for all $k$, and thus a.e.\ in $\R\times (0,\tau_1)$. 
If $\tau_1<T$, a finite iteration of this procedure proves uniqueness of the viscosity solution of $(CP)$.

Reassuming, we introduce a procedure, based on the barrier effect of spatial discontinuities, which indicates how the comparison result for semicontinuous sub- and supersolutions can be used to prove uniqueness of suitably defined viscosity solutions with discontinuous initial data. Although we have only done this for a particularly simple problem, preliminary calculations suggest that the procedure can be adapted to more general problems (to be addressed in future papers), namely the cases of initial data with infinitely many jump discontinuities and
Hamiltonians with linear growth and explicit $x$ and $t$ dependence. The latter case is particularly interesting since it includes equations for which uniqueness of discontinuous viscosity solutions fails, as mentioned in the beginning of the Introduction. In particular our approach seems to suggest a mathematical uniqueness criterium.
Many other problems concerning discontinuous solutions of first order HJ equations remain to be solved, in particular the multidimensional case.

Finally we observe that, setting  $v=u_x$ and $v_0=u_0'$, problem $(CP)$ is formally related to  the Cauchy problem for a scalar conservation law,
\begin{equation*}
\begin{cases}
v_t+ [H(v)]_x=0&\text{in }\R\times (0,T)\\
v=v_0&\text{in } \R\times\{0\}\,,
\end{cases}
\leqno{(CL)}
\end{equation*}
although it is not trivial to make the correspondence rigorous (\cite{BSTT10}; see also \cite{Ba} for a statement in this direction). If $v_0=u_0'$ is a signed Radon measure, it is possible to prove existence of suitably defined measure-valued entropy solutions of problem $(CL)$ (\cite{BSTT11}; see also \cite{BSTT13, BSTT12} for the case of positive initial measures). Remarkably, if the singular part $v_{0s}$ of $v_0$ (with respect to the Lebesgue measure) is a finite superposition of Dirac masses, the uniqueness of such solutions requires some  additional {\em compatibility conditions} to be satisfied near the support of $v_{0s}$. The singularities of $v_0$ in $(CL)$ have a barrier effect which corresponds to that produced by discontinuities of $u_0$ in $(CP)$, and well-posedness of $(CL)$ can be proven using {\em singular Dirichlet problems} which are the natural counterpart of singular Neumann problems for HJ equations (see \cite{BSTT12, BSTT11} for details).


\section{Semicontinuous envelopes}\label{prelim}
\setcounter{equation}{0}

Let $\chi_E$ denote the characteristic function of $E\subseteq\R$. For all $u\in\R$ we set
\begin{equation*}
[u]_\pm=\max\{\pm u,0\}, \quad{\rm sgn}_\pm(u)=\pm\chi_{\R_\pm}(u), \quad {\rm sgn}(u)={\rm sgn}_-(u)+{\rm sgn}_+(u)\,.
\end{equation*}

Let $D\subset\R^2$ be open, $z:D\mapsto\R$ a measurable function,  and $(x_0,t_0)\in \overline{D}$. We set 
$$
 \text{\rm ess}\!\!\!\!\!\!\!\!\!\limsup_{D\ni (x,t)\to (x_0,t_0)}\!\!z(x,t)\,:=\,\inf_{\delta>0}\, \left(\ \,{\rm ess}\!\!\!\!\!\!\!\!\! \!\!\!\!\sup_{(x,t)\in D\cap B_{\delta}(x_0,t_0)}\!\! z(x,t)\right)= \lim\limits_{\delta\to 0^+} \left(\ \, {\rm ess}\!\!\!\!\!\!\!\!\! \!\!\!\!\sup_{(x,t)\in D\cap B_{\delta}(x_0,t_0)}\!\! z(x,t)\right) ,
$$
$$
  \text{\rm ess}\!\!\!\!\!\!\!\!\!\liminf_{D\ni (x,t)\to (x_0,t_0)}\!\!z(x,t)\,:=\,\sup_{\delta>0} \,\left(\ \, {\rm ess}\!\!\!\!\!\!\!\!\! \!\!\!\!\inf_{(x,t)\in D\cap B_{\delta}(x_0,t_0)}\!\! z(x,t)\right)= \lim\limits_{\delta\to 0^+} \left(\ \, {\rm ess}\!\!\!\!\!\!\!\!\! \!\!\!\!\inf_{(x,t)\in D\cap B_{\delta}(x_0,t_0)}\!\! z(x,t)\right) ,
$$
where 
$$
B_r(x_0,t_0)=\{(x,t)\in \R^2\,|\, (x-x_0)^2+(t-t_0)^2< r^2\} \qquad(r>0)\,.
$$ 
If  $\text{\rm ess}\limsup_{D\ni (x,t)\to (x_0,t_0)}z(x,t)=   \text{\rm ess}\liminf_{D\ni (x,t)\to (x_0,t_0)}z(x,t)$, we also set 
$$
 \text{\rm ess}\!\!\!\!\!\!\!\!\!\!\!\!\lim_{D\ni (x,t)\to (x_0,t_0)}z(x,t)\,:= \, \text{\rm ess}\!\!\!\!\!\!\!\!\!\limsup_{D\ni (x,t)\to (x_0,t_0)}z(x,t)=   \text{\rm ess}\!\!\!\!\!\!\!\!\!\liminf_{D\ni (x,t)\to (x_0,t_0)}z(x,t)\,. 
$$ 
The quantities
$$
\text{\rm ess}\!\!\!\!\!\!\!\!\limsup_{D\ni (x,t)\to (x_0, t_0^{\pm})}z(x,t), \quad \text{\rm ess}\!\!\!\!\!\!\!\!\liminf_{D\ni (x,t)\to (x_0, t_0^{\pm})}z(x,t)
$$ 
are defined replacing $B_r(x_0,t_0)$ by $B_r(x_0,t_0) \cap \{(x,t)\in \R^2\,|\,t \gtrless t_0\}$. 
Similarly,
$$
\text{\rm ess}\!\!\!\!\!\!\!\!\limsup_{D\ni (x,t)\to (x_0^\pm,t_0)}z(x,t), \quad \text{\rm ess}\!\!\!\!\!\!\!\!\liminf_{D\ni (x,t)\to (x_0^\pm,t_0)}z(x,t)
$$ 
are defined replacing $B_r(x_0,t_0)$ by $B_r(x_0,t_0) \cap \{(x,t)\in \R^2\,|\,x \gtrless x_0\}$. 

\smallskip

Let $z\in L^\infty_{{\rm loc}}(\overline{D})$. By the {\it essential upper semicontinuous envelope} of $z$ we mean the function $z^*:\overline{D}\to \R$,
\begin{equation}\label{use}
z^*(x_0,t_0)\,:=\,\text{\rm ess}\!\!\!\!\!\!\!\!\!\limsup_{ D \ni (x,t)\to (x_0,t_0)}\!\! z(x,t) \quad\text{for any }(x_0,t_0)\in \overline{ D }\,.
\end{equation} 
The {\it essential lower semicontinuous envelope} of $z$ is the function  $z_*:\overline{ D }\to \R$, 
\begin{equation}\label{lse}
z_*(x_0,t_0)\,:=\,\text{\rm ess}\!\!\!\!\!\!\!\!\!\liminf_{ D \ni (x,t)\to (x_0,t_0)}\!\! z(x,t) \quad\text{for any }(x_0,t_0)\in \overline{ D }\,.
\end{equation}
We also set
\begin{equation}\label{limt}
z^*(x_0, t_0^{\pm})\,:=\,\text{\rm ess}\!\!\!\!\!\!\!\!\!\limsup_{ D \ni (x,t)\to (x_0,t_0^{\pm})}\!\!\!\! z(x,t)\,, \quad
z_*(x_0, t_0^{\pm})\,:=\,\text{\rm ess}\!\!\!\!\!\!\!\!\!\liminf_{ D \ni (x,t)\to (x_0,t_0^\pm)}\!\!\!\! z(x,t)\,,
\end{equation}
\begin{equation}\label{limx}
z^*(x_0^\pm,t_0)\,:=\,\text{\rm ess}\!\!\!\!\!\!\!\!\!\limsup_{ D \ni (x,t)\to (x_0^\pm,t_0)}\!\!\!\! z(x,t)\,, \quad
z_*(x_0^\pm,t_0)\,:=\,\text{\rm ess}\!\!\!\!\!\!\!\!\!\liminf_{ D \ni (x,t)\to (x_0^\pm,t_0)}\!\!\!\! z(x,t)\,.
\end{equation}
Observe that
\begin{equation}\label{lineq1}
z^*(x_0, t_0^{\pm})\le z^*(x_0,t_0)\,, \quad
z_*(x_0, t_0^{\pm})\ge z_*(x_0,t_0)\,,
\end{equation}
\begin{equation}\label{lineq2}
z^*(x_0^\pm,t_0)\le z^*(x_0,t_0)\,, \quad
z_*(x_0^\pm,t_0)\ge z_*(x_0,t_0)\,.
\end{equation}
Similar definitions hold for any measurable function $z:F\subseteq \R\to \R$.
For shortness, we shall say  ``upper (respectively lower) envelope'' instead of  ``essential upper (respectively lower) semicontinuous envelope''.  

It is easily checked that the upper (lower) envelope $z^*$ ($z_*$) is indeed upper (lower) 
semicontinuous in $\overline{ D }$, namely for any $(x_0,t_0)\in\overline{ D }$
\begin{equation}\label{uscu*}
\limsup\limits_{\overline{ D }\ni (x,t)\to (x_0,t_0)}\!z^*(x,t)\le z^*(x_0,t_0), \qquad
\liminf\limits_{\overline{ D }\ni (x,t)\to (x_0,t_0)}\! z_*(x,t)\ge  z_*(x_0,t_0).
\end{equation}
The inequalities in \eqref{uscu*} can be replaced by equalities. More generally we have: 

\begin{lem}\label{propis}
Let  $D\subseteq \R^2$ be  open and $z\in L^\infty_{{\rm loc}}(\overline{D})$. 
Then
\begin{equation}\label{el generale}
(z^*)^*=(z_*)^*=z^*,\quad 
(z^*)_*=(z_*)_*=z_*\quad\text{in }\overline D.
\end{equation}
\end{lem}

\begin{proof} 
Let $(x_0,t_0)\in \overline D$. We only prove that $(z^*)^*(x_0,t_0)=(z_*)^*(x_0,t_0)=z^*(x_0,t_0)$.
Since
$$
\begin{aligned}
(z_*)^*(x_0,t_0)&=\text{ess}\!\!\!\!\!\!\limsup_{ D \ni (x,t)\to (x_0,t_0)}z_*(x,t)\le (z^*)^*(x_0,t_0)\,=\\
&=\text{ess}\!\!\!\!\!\!\limsup_{ D \ni (x,t)\to (x_0,t_0)}z^*(x,t)
\le \limsup\limits_{\overline{ D }\ni (x,t)\to  (x_0,t_0)}z^*(x,t)\le z^*(x_0,t_0)
\end{aligned}
$$
(see \eqref{use}, \eqref{lse} and \eqref{uscu*}), it is enough to show that 
\begin{equation}\label{<limsup}
z^*(x_0,t_0)\le \text{ess}\!\!\!\!\!\!\limsup_{\overline{ D }\ni (x,t)\to (x_0,t_0)} z_*(x,t)\,.
\end{equation}
For every $\ep>0$ there exists $r_{\ep}>0$ such that 
$${\rm ess}\!\!\! \!\!\!\!\!\!\!\!\!\!\!\sup_{(x,t)\in B_{r}(x_0,t_0)\cap  D }\!\! z(x,t) > z^*(x_0,t_0)-\tfrac12\ep\quad\mbox{for all}\ \,r \leq r_{\ep}\,.$$
Therefore, for every such $r$ there exists $\mathcal{B}_{r,\ep}\subseteq B_{r}(x_0,t_0)\cap  D $, $|\mathcal{B}_{r,\ep}|>0$, such that 
$$
z(x,t)\geq {\rm ess}\!\!\! \!\!\!\!\!\!\!\!\!\!\!\sup_{(x,t)\in B_{r}(x_0,t_0)\cap  D }\!\! z(x,t) -\tfrac12\ep > z^*(x_0,t_0)-\ep\quad\ \mbox{for a.e.~$(x,t)\in \mathcal{B}_{r,\ep}$}\,.$$
Setting $\mathcal{B}_\ep=\bigcup_{r\leq r_{\ep}}\mathcal{B}_{r,\ep}$, it follows that 
$|\,\mathcal{B}_\ep\cap B_r(x_0,t_0)|>0$ for every $r>0$ 
(hence $|\mathcal{B}_{\ep}|>0$), and
$$
z>z^*(x_0,t_0)-\ep \quad  \text{a.e.~in }\mathcal{B}_\ep\,.
$$

Let  $(x,t)\in \mathcal{B}_\ep$ satisfy $|B_{\delta}(x,t)\cap \mathcal{B}_\ep|>0$ for every $\delta>0$; this choice is possible up to a null set $\mathcal{N}_{\ep}\subset\mathcal{B}_{\ep}$, since $|\mathcal{B}_{\ep}|>0$ and almost every $(x,t)\in \mathcal{B}_\ep$ is a Lebesgue point of $f(x,t)=\chi_{\mathcal{B}_{\ep}}(x,t)$ ($e.g.$, see \cite[Subsection 1.7.1]{EG}). Then we have
$$
z>z^*(x_0,t_0)-\ep \quad\text{a.e.~in }\mathcal{B}_\ep\cap B_{\delta}(x,t), \qquad (|\,\mathcal{B}_\ep\cap B_{\delta}(x,t)|>0\quad \mbox{for all}\ \delta>0)\,,
$$
whence $z_*(x,t)\geq z^*(x_0,t_0)-\ep$ for a.e.~$(x,t)\in \mathcal{B}_\ep$. Since $|\,\mathcal{B}_\ep\cap B_r(x_0,t_0)\,|>0$ for all $r>0$, this implies that 
$$
\text{ess}\!\!\limsup_{\overline{ D }\ni (x,t)\to (x_0,t_0)}z_*(x,t)\,\ge\,z^*(x_0,t_0)-\ep\,.
$$
and \eqref{<limsup} follows from the arbitrariness of $\ep$.
\end{proof}
Similar results  hold if $D\subseteq \R$ is open and $z\in L^\infty_{{\rm loc}}( \overline{D} )$. 

\smallskip

\begin{remark}\label{restri}
Since the definition of $z^*,z_*$ depends on the domain of definition of $z$, different restrictions of $z$ can have different upper and lower envelopes. In fact, let $ D _1\subseteq  D $ be an open set, and let $z_1= z\lefthalfcup  D _1$ be the restriction of $z$ to $ D _1$. 
If $(x_0,t_0)\in\partial  D _1$, then:
$$
(z_1)^*(x_0,t_0)=\text{\rm ess}\!\!\!\!\!\!\!\!\!\limsup_{ D _1\ni (x,t)\to (x_0,t_0)}\!\! z(x,t) \leq \text{\rm ess}\!\!\!\!\!\!\!\!\!\limsup_{ D \ni (x,t)\to (x_0,t_0)}\!\! z(x,t) =z^*(x_0,t_0)\,, 
$$
$$
(z_1)_*(x_0,t_0)=\text{\rm ess}\!\!\!\!\!\!\!\!\!\liminf_{ D _1\ni (x,t)\to (x_0,t_0)}\!\! z(x,t) \geq \text{\rm ess}\!\!\!\!\!\!\!\!\!\liminf_{ D \ni (x,t)\to (x_0,t_0)}\!\! z(x,t) = z_*(x_0,t_0)
$$
(of course, if $(x_0,t_0)\in\mathring{ D }_1$ these inequalities become equalities). Observe that \eqref{lineq1}-\eqref{lineq2} are a particular case of these inequalities. If $D=\O\times(0,T)$ and $t_0=0$, inequalities \eqref{lineq1} from above become 
equalities, namely
\begin{equation}\label{lineq0}
z^*(x_0, 0^+)= z^*(x_0,0)\,, \quad
z_*(x_0, 0^+)= z_*(x_0,0)\quad\text{for any $x_0\in\overline{\O}$}\,.
\end{equation}

More precisely, let $ D _1\subseteq  D $ be an open set, and let $z_1= z\lefthalfcup  D _1$. 
Set 
\begin{equation}\label{tilq}
\tilde{ D }_1\,:=\,\{(x,t)\in \overline{ D }_1\,|\, \exists\, \d_1>0\; \text{such that} \;B_{\d_1}(x,t)\cap  D \,\subseteq\, D _1\}\,.
\end{equation}
Then for any $(x_0,t_0)\in \tilde{ D }_1$ there holds
\begin{equation}\label{tilque}
z^*(x_0,t_0)=(z_1)^*(x_0,t_0),\qquad z_*(x_0,t_0)=(z_1)_*(x_0,t_0)\,.
\end{equation}
In fact, if $(x_0,t_0)\in \tilde{ D }_1$ then $ D \cap B_{\d_1}(x_0,t_0) = D _1\cap B_{\d_1}(x_0,t_0)$ and,
by \eqref{use},
$$ z^*\!(x_0,t_0)
\!=\! \lim\limits_{\delta\to 0^+} \!\left( {\rm ess}\!\!\!\!\!\! \!\!\!\!\sup_{(x,t)\in  D \cap B_{\delta}(x_0,t_0)}\!\!\!\!\! z(x,t)\!\right)\!=\!\lim\limits_{\delta\to 0^+} \!\left( {\rm ess}\!\!\!\!\!\! \!\!\!\!\sup_{(x,t)\in  D _1\cap B_{\delta}(x_0,t_0)}\!\! \!\!\!z_1(x,t)\!\right)\!=\!
(z_1)^*\!(x_0,t_0).
$$

For further reference we consider some specific cases of $D_1\subseteq D=(a,b)\times(0,T)$, 
with $-\infty<a<b<\infty$. The first example concerns a trapezoidal domain.
\begin{itemize}
\item[$(a)$] Let
$
D _1= \{(x,t)\,|\,x\in(a,b-\|H'\|_{\infty}t),\,t\in (0,\tau)\} 
$  
for $\tau\in (0,\tau_1]$, where 
\begin{equation}\label{det1}
\t_1=\min\left\{\tfrac12(b-a)/\|H'\|_\infty,\, T\right\}\,.
\end{equation}
Then 
$\tilde{ D }_1= \{(x,t)\,|\,x\in[a,b-\|H'\|_{\infty}t),\,t\in [0,\tau)\}$ if $\tau\in (0,\tau_1)$, and $\tilde{ D }_1= \{(x,t)\,|\,x\in[a,b-\|H'\|_{\infty}t),\,t\in [0,T]\}$ if $\tau=\tau_1=T$; \smallskip

\item[$(b)$] let $ D _1=Q_{\t,T}\,:=\,\{(x,t)\in D\,|\, t\in(\t,T)\}$ for any 
$\tau\in (0,T)$. Then
$\tilde{ D }_1= \{(x,t)\,|\,x\in[a,b],\,t\in (\t,T]\}$; \smallskip

\item[$(c)$] let $ D _1=\{(x,t)\in D\,|\, x\in (a,c),\, t\in(0,T)\}$ for any $c\in(a,b)$. Then
$\tilde{ D }_1= \{(x,t)\,|\,x\in[a,c),\,t\in [0,T]\}$.
\end{itemize}
\end{remark}


\section{Definitions and results}\label{res}
\setcounter{equation}{0}

As we have explained in the Introduction, the proof of the main uniqueness result requires the introduction of {\it singular Neumann problems}.
Given $a,b\in \R$, $a<b$, we consider the problems
$$
\begin{aligned}
&(N_{a,\pm}) \
\begin{cases}
u_t+ H(u_x)=0&\text{in }(a,\infty)\times (0,T)\\
u_x(a,t)=\pm\infty&\text{for } t\in(0,T),
\end{cases}
\\
&(N^{b,\pm}) \
\begin{cases}
u_t+ H(u_x)=0&\text{in }(-\infty,b)\times (0,T)\\
u_x(b,t)=\pm\infty&\text{for } t\in(0,T),
\end{cases}
\\
&(N_{a,\pm}^{b,\pm}) \
\begin{cases}
u_t+ H(u_x)=0&\text{in }(a,b)\times (0,T)\\
u_x(a,t)=\pm\infty&\text{for } t\in(0,T)\\
u_x(b,t)=\pm\infty&\text{for } t\in(0,T).
\end{cases}
\\
\end{aligned}
$$
In addition we consider equation \eqref{hje} in the following trapezoidal domains: 
\begin{equation}\label{trap bis}
\begin{cases}
A\,:=\,\left\{(x,t)\,|\, x\!\in\!(a,d-\|H'\|_{\infty}t),\,t\in (0,T)\right\} &\text{with $d-a\ge\|H'\|_{\infty}T $}\\
B\,:=\,\left\{(x,t)\,|\, x\!\in\!(c+\|H'\|_{\infty}t,b),\,t\in (0,T)\right\}& \text{with $\|H'\|_{\infty}T\le b-c$}\\
C\,:=\,\left\{(x,t)\,|\, x\!\in\!(c\!+\!\|H'\|_{\infty}t,d\!-\!\|H'\|_{\infty}t),t\in (0,T)\right\} &\text{with $d-c\ge2\|H'\|_{\infty}T$}.
\end{cases}
\end{equation} 
Due to their slope, no boundary conditions are required on the oblique sides:
$$
\begin{aligned}
&({\rm T}_{a,\pm}) \
\begin{cases}
u_t+ H(u_x)=0&\text{in }A\\ 
u_x(a,t)=\pm\infty&\text{for } t\in(0,T),
\end{cases}
\\&
({\rm T}^{b,\pm}) \
\begin{cases}
u_t+ H(u_x)=0&\text{in }B\\
u_x(b,t)=\pm\infty&\text{for } t\in(0,T),
\end{cases}
\\
&({\rm T}_C) \quad 
u_t+ H(u_x)=0\quad\ \text{in }C.
\\
\end{aligned}
$$

All these problems require  an initial condition $u_0$, which is assumed to be locally bounded in the main existence and regularity results (see Proposition \ref{prolip} and Theorem \ref{exi} below) and, in addition, 
piecewise continuous in the uniqueness Theorem \ref{unique}.
More precisely, a function $f\in L^{\infty}(\O)$ defined in a bounded interval $\O=(a,b)$ is said to be {\em piecewise continuous in $\O$} if
\smallskip

\noindent - $\O=\bigcup_{j=1}^{p+1} I_j$ $(p\in \N)$ with $I_1=(a,x_1)$, $I_j=(x_{j-1},x_j)$ if $2\le j\le p$, $I_{p+1}=(x_p,b)$; 
\smallskip

\noindent - the restriction $f\lefthalfcup I_j$ has a representative $f_j\in C(\overline{I}_j)$ 
for $1\le j\le p+1$, $f_j(x_j)\neq f_{j+1}(x_j)$ for $1\le j\le p$. 
\smallskip

\noindent Observe that a function $f\in C(\overline{\O})$ is piecewise continuous (the case $p=0$).
If instead $\Omega$ is an unbounded interval, a function $f\in L^\infty_{{\rm loc}}(\overline{\O})$ is piecewise continuous in $\O$ 
if it is piecewise continuous in every bounded interval $(a_0,b_0 )\subset \O$. 
\smallskip

\subsection{Definitions}\label{def}
The following definitions  are used throughout the paper.
\begin{definition}\label{visceq}
Let $D\subseteq \R^2$ open and $u\in L^\infty_{{\rm loc}}(\overline{D})$.  

\noindent $(i)$ $u$ is a {\em viscosity subsolution of equation \eqref{hje} in $D$} if  for all $\varphi\in C^1(D)$ the following condition holds: 
if $(x,t)\!\in\! D$ is a local maximum point of $u^*\!-\varphi$ in $D$, then
\begin{equation}\label{ad0-}
\varphi_t(x,t)+H(\varphi_x(x,t))\le 0\,. 
\end{equation} 

\noindent $(ii)$ $u$ is a {\em viscosity supersolution of equation \eqref{hje} in $D$} if for all $\varphi\in C^1(D)$ the following condition holds: 
if $(x,t)\in D$ is a local minimum point of $u_*-\varphi$ in $D$, then
\begin{equation}\label{ad0+}
\varphi_t(x,t)+H(\varphi_x(x,t))\ge 0\,. 
\end{equation} 
\end{definition}
If $D\subseteq \R^2$ is bounded, viscosity sub and supersolutions of \eqref{hje} in $D$ belong to $L^{\infty}(D)$.
\smallskip

We first define the solution concept for problem $(N_{a,\pm}^{b,\pm})$, and then explain how to extend it to other problems.    

\begin{definition}\label{defsubsupervi}
Let $-\infty< a < b<\infty$, $Q=(a,b)\times (0,T)$ and $\hat Q=[a,b]\times (0,T]$. 
\smallskip

\noindent $(i)$ 
A viscosity subsolution $u$ of \eqref{hje} in $Q$ is a {\it viscosity subsolution} of:
\smallskip

- problem $(N_{a,-}^{b,+})$;
\smallskip

- problem $(N_{a,+}^{b,+})$, if for all $\varphi\in C^1(\hat{Q})$ 
 \begin{equation}\label{ad1-}
u^*-\varphi \text{ has a local maximum at }(a,t)\in \hat Q\ \Rightarrow\ \varphi_t(a,t)+H(\varphi_x(a,t))\le 0; 
\end{equation} 

- problem $(N_{a,-}^{b,-})$, if for all $\varphi\in C^1(\hat{Q})$ 
\begin{equation}\label{ad2-}
u^*-\varphi \text{ has a local maximum at }(b,t)\in \hat Q\ \Rightarrow\ \varphi_t(b,t)+H(\varphi_x(b,t))\le 0; 
\end{equation} 

- problem $(N_{a,+}^{b,-})$, if for all $\varphi\in C^1(\hat{Q})$ both  \eqref{ad1-} and \eqref{ad2-} are satisfied.
\smallskip

\noindent $(ii)$
A viscosity supersolution $u$ of \eqref{hje} in $Q$ is a {\it viscosity supersolution} of:
\smallskip

- problem $(N_{a,+}^{b,-})$;
\smallskip

- problem $(N_{a,+}^{b,+})$, if for all $\varphi\in C^1(\hat{Q})$ 
 \begin{equation}\label{ad1+}
u^*-\varphi \text{ has a local minimum at }(b,t)\in \hat Q\ \Rightarrow\ \varphi_t(b,t)+H(\varphi_x(b,t))\ge 0; 
\end{equation} 

- problem $(N_{a,-}^{b,-})$, if for all $\varphi\in C^1(\hat{Q})$ 
\begin{equation}\label{ad2+}
u^*-\varphi \text{ has a local minimum at }(a,t)\in \hat Q\ \Rightarrow\ \varphi_t(a,t)+H(\varphi_x(a,t))\ge 0; 
\end{equation} 

- problem $(N_{a,-}^{b,+})$, if for all $\varphi\in C^1(\hat{Q})$ both  \eqref{ad1+} and \eqref{ad2+} are satisfied.
\smallskip

\noindent $(iii)$ A function $u$ is a {\em viscosity solution} of problem $(N_{a,\pm}^{b,\pm})$ 
if it is both a viscosity subsolution and a viscosity supersolution of $(N_{a,\pm}^{b,\pm})$.

\smallskip

\noindent $(iv)$ Let $u_0\in L^\infty(a,b)$. A {\em viscosity solution of $(N_{a,\pm}^{b,\pm})$ with initial condition $u_0$} is a viscosity solution of $(N_{a,\pm}^{b,\pm})$ such that 
\begin{equation}\label{cic}
u^*(\cdot,0)=(u_0)^*\,, \quad u_*(\cdot,0)= (u_0)_* \;\;\text{ in $[a,b]$}\,. 
\end{equation}
\end{definition}
\medskip

Observe that Definition \ref{defsubsupervi} makes sense for any $H\in C(\R)$.

It is straightforward to modify Definition \ref{defsubsupervi} in case of the other problems we mentioned before.
For example, a viscosity subsolution of problem (${\rm T}_{a,-}$) is a viscosity subsolution of \eqref{hje} in the set $A$, a viscosity subsolution of problem (${\rm T}_{a,+}$) is a viscosity subsolution of \eqref{hje} in $A$ which satisfies \eqref{ad1-} for all $\varphi\in C^1(\hat A)$ (where $\hat A=\overline A\cap (\R\times (0,T])$), and a viscosity subsolution of problem (${\rm T}_C$) is a viscosity subsolution of \eqref{hje} in the set $C$. 

If $u_0$ is defined in an unbounded interval, in point $(iv)$ we only require boundedness of $u_0$ in all bounded subintervals. In particular: 
\begin{definition}
For every $u_0\in L_{\rm loc}^\infty(\R)$, $u$ is a {\it viscosity solution of the Cauchy problem $(CP)$} if it is a viscosity sub- and supersolution of \eqref{hje} in $\R\times (0,T)$ and satisfies \eqref{cic} in $\R$.
\end{definition}
\smallskip

\begin{remark} \label{remres}
Let $Q=(a,b)\times(0,T)$ with $-\infty<a<b<\infty$, and $Q_1=(a,b)\times (\tau_1,\tau_2)$, for $0\leq \tau_1<\tau_2\leq T$. 
It is easily seen that if $u$ is a viscosity subsolution of $(N_{a,\pm}^{b,\pm})$  in $Q$, 
its restriction $u_1= u\lefthalfcup Q_1$  is a viscosity subsolution of $(N_{a,\pm}^{b,\pm})$  in $Q_1$.
In fact, if $(x_0,t_0)$ is a local maximum point of $(u_1)^*-\varphi$ in $\hat{Q}_1$ and $t_0\in (\tau_1,\tau_2)$, by \eqref{tilque} it is also a local maximum point of 
$u^*-\varphi$ in $\hat{Q}$, whereas in the case $t_0=\tau_2$ it suffices to argue as in \cite[Section 5.2]{Ev}. 
 
Similar remarks hold for viscosity subsolutions of equation \eqref{hje} in $Q$
or of any of the other problems mentioned before, as well as for viscosity supersolutions.
\end{remark}


\subsection{Main results}\label{exist} 
In Section \ref{compre} we prove the basic comparison result for viscosity sub- and supersolutions: 
\begin{theorem}\label{comp} Let $(H_1)$ hold. Let $u$ be a viscosity subsolution and $v$ a viscosity supersolution of problem $(N_{a,\pm}^{b,\pm})$.
Then
\begin{equation}\label{cfrt}
\max_{[a,b]\times [0,t]} [u^*-v_*]_+\, \le  \, \max_{[a,b]}\, \left[u^*(\cdot,0)-v_*(\cdot,0)\right]_+\quad\text{for all }t\in (0,T].
\end{equation}
Similar results are valid for problems $({\rm T}_{a,\pm})$, $({\rm T}^{b,\pm})$ and $({\rm T}_C)$: 
\begin{equation}\label{cfrtbis}
\begin{cases}
\max_{\overline A}\,[u^*-v*]_+\, \le  \, \max_{[a,d] }\, \left[u^*(\cdot,0)-v_*(\cdot,0)\right]_+ &\text{for problem $({\rm T}_{a,\pm})$}\\
\max_{\overline B}\,[u^*-v*]_+\, \le  \, \max_{[c,b] }\, \left[u^*(\cdot,0)-v_*(\cdot,0)\right]_+ &\text{for problem $({\rm T}^{b,\pm})$}\\
\max_{\overline C}\,[u^*-v*]_+\, \le  \, \max_{[c,d] }\, \left[u^*(\cdot,0)-v_*(\cdot,0)\right]_+ &\text{for problem $({\rm T}_C)$.}\\
\end{cases}
\end{equation}
\end{theorem}
Theorem \ref{comp} will be proven by a method of doubling variables adapted from \cite{PLL}, where only the case of homogeneous Neumann boundary conditions was considered (see also \cite{Ba2} for more general Neumann boundary conditions). It implies uniqueness and continuity of viscosity solutions of, for example, problem $(CP)$ if $u_0$ is continuous in $\R$ (since $u^*\equiv u_*$ in 
$\R\times [0,T]$). If instead $u_0$ is not continuous at a point $x=c$,
the following regularity results give information about the temporal evolution of $u^*-u_*$ at $c$. These result, proved in Section~\ref{Sec reg}, are crucially based on the boundedness of $H$. 

\begin{prop}\label{prolip} Let $(H_1)$ hold and let $\O=(a,b)$, with $-\infty\le a<c<b\le \infty$, let
\begin{equation}\label{Def K,k}
K:=\sup_{s\in \R}\,[-H(s)]\,, \qquad k:=\inf_{s\in \R}\,[-H(s)]\,,
\end{equation}
and let $u_0\in L^\infty_{\rm loc}(\overline \O)$ be such that the one-sided essential limits $u_0(c^\pm)$ exist.
\smallskip

\noindent $(i)$
Let $u$ and $v$ be a viscosity sub- and supersolution of equation \eqref{hje} in $(a,b)\times (0,T)$. 
 Then for any $t_1,t_2\in [0,T]$, $t_1\ne t_2$, and $x\in\overline {\O}$
\begin{equation}\label{sih}
\frac{u^*(x,t_1)-u^*(x,t_2)}{t_1-t_2}\le  K, \qquad
\frac{v_*(x,t_1)-v_*(x,t_2)}{t_1-t_2}\ge  k\,.
\end{equation}
\noindent $(ii)$
Let $u$ be a viscosity solution of equation \eqref{hje} in $(a,b)\times (0,T)$ with initial function $u_0$. Let $u_0(c^+)\neq u_0(c^-)$ and set
$$
\underline t\,:=\left\{\begin{array}{ll}
\frac{|u_0(c^+)-u_0(c^-)|}{K-k}&\mbox{if}\ \,K>k\,,\smallskip\\
T&\mbox{if}\ \,K=k\,.
\end{array}\right.
$$
Then, for all $t\in(0, \min\{\underline t,T\})$,
\begin{equation}\label {wt 1}
\begin{cases}
u^*(c,t)=u^*(c^+,t)>u^*(c^-,t)\\
u_*(c,t)=u_*(c^-,t)<u_*(c^+,t)
\end{cases}
\quad \text{if } u_0(c^+)> u_0(c^-),
\end{equation}
\begin{equation}\label {wt 2}
\begin{cases}
u^*(c,t)=u^*(c^-,t)>u^*(c^+,t)\\
u_*(c,t)=u_*(c^+,t)<u_*(c^-,t)
\end{cases}
\quad \text{if } u_0(c^+)< u_0(c^-).
\end{equation}
If  \eqref {wt 1}  or \eqref {wt 2} holds for $0<t<t_c$ for some $t_c\in [\min\{\underline t,T\},T]$, then for every $0\leq t_0<t_1<t_c$ we have
\begin{subequations}\label{desa1}
\begin{equation}\label{desa+}
u^*(c,t_1)-u_*(c,t_1)\leq 
u^*(c,t_0)-u_*(c,t_0) - A_+(t_1-t_0)\quad\text{if $u_0(c^+)> u_0(c^-)$\,,}
\end{equation}
respectively
\begin{equation}\label{desa-}
u^*(c,t_1)-u_*(c,t_1)\leq 
u^*(c,t_0)-u_*(c,t_0) - A_-(t_1-t_0)\quad\text{if $u_0(c^+)< u_0(c^-)$\,,}
\end{equation}
\end{subequations}
where
\begin{equation}\label{A+-}
A_\pm\,:=\,\limsup\limits_{s\to \pm\infty}H(s)-\liminf\limits_{s\to \pm \infty}H(s)\,.  
\end{equation} 
\end{prop}

Part $(i)$, a sort of Lipschitz continuity with respect to $t$, might be easily guessed from the boundedness of $H$. Part $(ii)$
quantifies the fact  that initial jump discontinuities do not disappear instantaneously; in addition, if $A_\pm>0$ it will be seen below that for initial data as in $(H_2)$
spatial jumps may disappear in the interval $(0,T)$ (see Theorem \ref{unique} and Remark \ref{resa}; see also \cite{BSTT10}). 

In Section \ref{perron} we use Perron's method to prove the existence of viscosity solutions:

\begin{theorem}\label{exi} 
If $H$ satisfies $(H_1)$ and $u_0\in  L^\infty_{{\rm loc}}(\R)$, problem $(CP)$ has a viscosity solution.
Similarly, if $-\infty<a<b<\infty$, there exists a viscosity solution of: 

- problem $(N_{a,\pm})$ with initial function $u_0\in L^{\infty}_{\rm loc}([a,\infty))\,;$

-  problem $(N^{b,\pm})$ with initial function $u_0\in L^{\infty}_{\rm loc}(-\infty,b])\,;$

- problem $(N_{a,\pm}^{b,\pm})$ with initial function $u_0\in L^{\infty}((a,b))\,$.
\end{theorem} 

Theorem \ref{comp} will be used to prove uniqueness of discontinuous viscosity solutions of problems $(CP)$, $(N_{a,\pm})$, $(N^{b,\pm})$ and $(N_{a,\pm}^{b,\pm})$ with initial data which have a finite number of jump discontinuities. Below we indicate these problems generically by 
$$
\text{``problem $(N)$ with initial condition $u_0$ satisfying $(H_2)$ in $\O=(a,b)$''.}
$$ 
So $-\infty\le a < b\le\infty$, $a=-\infty$ in case of problems $(CP)$ and $(N^{b,\pm})$, $b=\infty$ in case of problems $(CP)$ and $(N_{a,\pm})$, and in condition $(H_2)$ the set $\R$ is replaced by $\O$.

We  denote by $x_j$ the points where $u_0$ is discontinuous and set $Q=(a,b)\times (0,T)$ and $Q_j=I_j\times (0,T)$ for $j=1,\dots,p+1$, where
$$
I_1\,:=\,(a,x_1), \qquad I_j\,:=\,(x_{j-1},x_j) \quad \text{if }j=2,\dots,p, \qquad I_{p+1}\,:=\,(x_p,b)\,.
$$ 
The main statement of the following result is part $(i)$, the uniqueness of viscosity solutions. 

\begin{theorem}\label{unique} Let $H$ satisfy $(H_1)$, and let $K$, $k$ be defined by \eqref{Def K,k}. Let $u_0$ satisfy $(H_2)$ in $(a,b)$ $(-\infty\le a<b\le\infty)$.

\noindent $(i)$ If $u$ and $v$ are viscosity solutions of problem $(N)$ with initial condition $u_0$, then $u=v$ a.e.~in $Q$.

\noindent $(ii)$ If $u$ is the viscosity solution of problem $(N)$ with initial condition $u_0$, then:
\begin{itemize}
\item[$(a)$]  for all $j=1,\dots,p+1$ the restriction $u_j=u\lefthalfcup Q_j$ has a continuous representative $\tilde{u}_j$ in $\overline{Q}_j$; 
\item[$(b)$] for all $j=1,\dots,p$ there exists a unique $\t_j\in(0,T]$ such that 
$$
\tilde{u}_j(x_j,t)\neq \tilde{u}_{j+1}(x_j,t) \;\;\Leftrightarrow\;\; t\in[0,\t_j)\,; 
$$
\item[$(c)$] for all $j=1,\dots,p+1$ the representative $\tilde{u}_j$ is Lipschitz continuous with respect to $t$ in $\overline{Q}_j$:  
\begin{equation}\label{inelip}
k \leq \frac{\tilde{u}_j(x,t_1)-\tilde{u}_j(x,t_2)}{t_1-t_2}\le  K\quad\text{for all $t_1,t_2\in[0,T ]$, $t_1\neq t_2$, and $x\in \overline{I}_j$.}
\end{equation}
\item[$(d)$] for all $j=1,\dots,p$ and $0\le t_0<t_1<\t_j$,  
$$
|\tilde u_{j+1}(x_j,t_1)-\tilde u_j(x_j,t_1)|\leq |\tilde u_{j+1}(x_j,t_0)-\tilde u_j(x_j,t_0)| 
-A_+ (t_1-t_0) \quad \mbox{if $u_0(x_j^+)>u_0(x_j^-)$},
$$
$$
|\tilde u_{j+1}(x_j,t_1)-\tilde u_j(x_j,t_1)|\leq |\tilde u_{j+1}(x_j,t_0)-\tilde u_j(x_j,t_0)| 
-A_- (t_1-t_0) \quad \mbox{if $u_0(x_j^-)>u_0(x_j^+)$},
$$
where $A_\pm$ is defined by \eqref{A+-} with $c$ replaced by $x_j$.
\end{itemize}
\end{theorem}

\begin{remark}\label{resa}
Let $(H_1)$-$(H_2)$ be satisfied and let $u$ be the viscosity solution of problem $(N)$ with initial function $u_0$. 
By claim $(a)$ in Theorem \ref{unique}$(ii)$,  
\begin{equation}
u^*(x,t)=u_*(x,t)\quad\mbox{for all $t\in [0,T]$, if $u_0$ is continuous at $x$.}
\end{equation}

Conversely, if $u_0$ has a jump discontinuity at $x_j$, then from claim $(b)$ in Theorem \ref{unique}$(ii)$ it follows that for every $t\in (0,T]$
\begin{equation}\label{eq jump}
u^*(x_j,t)-u_*(x_j,t)=|\tilde{u}_j(x_j,t)-\tilde{u}_{j+1}(x_j,t)|\,,
\end{equation}
and there exists $\t_j\in(0,T]$ such that
\begin{equation}\label{sal man}
\begin{cases}
u^*(x_j,t)-u_*(x_j,t)\neq 0&\mbox{for}\ \,t\in [0,\t_j)\\
u^*(x_j,t)-u_*(x_j,t)=0&\mbox{for}\ \,t\in [\t_j,T ]\quad\mbox{if}\ \,\tau_j<T\,.
\end{cases}
\end{equation}
More precisely, if $u_0$ has a jump discontinuity at $x_j$, then initial jumps cannot change sign: for all $t\in [0,\tau_j)$ there holds
\begin{equation}\label{salt0}
\left\{\begin{array}{ll}
u^*(x_j,t)=u^*(x_j^+,t)=\tilde{u}_{j+1}(x_j,t)>\tilde{u}_j(x_j,t)= u^*(x_j^-,t)\\
u_*(x_j,t)=u_*(x_j^-,t)=\tilde{u}_{j}(x_j,t)<\tilde{u}_{j+1}(x_j,t)=u_*(x_j^+,t)
\end{array}\right.
\ \ \mbox{if $u_0(x_j^+)>u_0(x_j^-)$}\,,
\end{equation}
\begin{equation}\label{salt0 bis}
\left\{\begin{array}{ll}
u^*(x_j,t)=u^*(x_j^-,t)=\tilde{u}_j(x_j,t)>\tilde{u}_{j+1}(x_j,t)= u^*(x_j^+,t)\\
u_*(x_j,t)=u_*(x_j^+,t)=\tilde{u}_{j+1}(x_j,t)<\tilde{u}_j(x_j,t)=u_*(x_j^-,t)
\end{array}\right.
\ \ \mbox{if $u_0(x_j^+)<u_0(x_j^-)$}
\end{equation}
(see also Proposition \ref{prolip} and Lemma \ref{ledi} below).
\end{remark}

An almost immediate consequence of Theorems \ref{exi} and \ref{comp} is a comparison principle for piecewise continuous viscosity solutions:

\begin{corollary}\label{thcs}
Let $H$ satisfy $(H_1)$, and let $u_0$ and $v_0$ satisfy $(H_2)$ in $\O=(a,b)$.
If $u$ and $v$ are viscosity solutions of problem $(N)$ 
with initial data $u_0\leq v_0$ a.e.~in $\Omega$, then $u\leq v$ a.e.~in $Q$.
\end{corollary}


\section{Comparison between viscosity sub- and supersolutions}\label{compre}
\setcounter{equation}{0}

In this Section we prove Theorem \ref{comp}. To do so, we need two preliminary results of independent interest.

\begin{prop}\label{p2} Let $H\in C(\R)$. 
Let $Q_1=(c,d)\times(t_1, t_2)$, $\hat{Q}_1=[c,d]\times(t_1, t_2]$  with $-\infty< c<d< \infty$, $0\le t_1< t_2\le T$. 

\noindent $(i)$ Let $u$ be a viscosity subsolution of equation \eqref{hje} in $Q_1$, and let $\varphi\in C^1(\hat{Q}_1)$. 
\begin{itemize}
\item[$(a)$] 
Let $(c,t_0)\in \hat{Q}_1$ be a local maximum point of $u^*-\varphi$ in $\hat{Q}_1$. Then
\begin{equation}\label{adter}
\varphi_t(c,t_0)+\inf\{H(\xi)\,|\,\xi\le \varphi_x(c^+,t_0)\}\le 0\,. 
\end{equation} 
\item[$(b)$] Let $(d,t_0)\in \hat{Q}_1$ be a local maximum point of $u^*-\varphi$ in $\hat{Q}_1$. Then
\begin{equation}\label{adterbis}
\varphi_t(d,t_0)+\inf\{H(\xi)\,|\,\xi\ge \varphi_x(d^-,t_0)\}\le 0\,. 
\end{equation} 
\end{itemize}

\noindent $(ii)$  Let $u$ be a viscosity supersolution of equation \eqref{hje} in $Q_1$, and let $\varphi\in C^1(\hat{Q}_1)$.
\begin{itemize}
\item[$(a')$] Let $(c,t_0)\in \hat{Q}_1$ be a local minimum point of $u_*-\varphi$ in $\hat{Q}_1$. Then
\begin{equation}\label{adqua}
\varphi_t(c,t_0)+\sup\{H(\xi)\,|\,\xi\ge \varphi_x(c^+,t_0)\}\ge 0\,. \end{equation} 
\item[$(b')$] Let $(d,t_0)\in \hat{Q}_1$ be a local minimum point of $u_*-\varphi$ in $\hat{Q}_1$. Then
\begin{equation}\label{aqb}
\varphi_t(d,t_0)+\sup\{H(\xi)\,|\,\xi\le \varphi_x(d^-,t_0)\}\ge 0\,. 
\end{equation} 
\end{itemize}
\end{prop}

\begin{proof} We only prove \eqref{adter}. The proofs of \eqref{adterbis}-\eqref{aqb} are similar.
 
Let $(c,t_0)$ be a local maximum point of $u^*-\varphi$ in $\hat{Q}_1$. Then 
\begin{equation*}
\limsup_{\overline Q_1\ni(y,\t)\to (c,t_0)}
\frac{u^*(y,\t)-u^*(c,t_0)-\varphi_x(c^+,t_0)(y-c)-\varphi_t(c,t_0)(\t-t_0)}{[(y-c)^2+ (\t-t_0)^2]^{1/2}}\le0,
\end{equation*} 
whence 
\begin{equation}\label{x2}
\limsup_{\overline Q_1\ni(y,\t)\to (c,t_0)}
\frac{u^*(y,\t)-u^*(c,t_0)-\xi(y-c)-\varphi_t(c,t_0)(\t-t_0)}{[(y-c)^2+ (\t-t_0)^2]^{1/2}}\le0
\end{equation} 
for all $\xi\ge\varphi_x(c^+,t_0)$ and
$$
\bar\xi=\inf\{\xi\in \R\,|\, \text{\eqref{x2} holds}\}\le\varphi_x(c^+,t_0)\,.
$$

First we consider the case that $\bar\xi>-\infty$. Then \eqref{x2} holds with $\xi=\bar\xi$ and
there exists $\psi\in C^1(\hat{Q}_1)$ such that $\psi_x(c,t_0)=\bar\xi$, $\psi_t(c,t_0)=\varphi_t(c,t_0)$, and 
$u^*-\psi$ has a strict maximum at $(c,t_0)$ ($e.g.$, see \cite[Proposition 2.6]{K}). 
For all sufficiently small $\delta>0$ the function  $u^*-\psi+\delta (x-c)$ has a maximum at some point $(x_\delta,t_\delta)\in (c,d)\times (t_1,t_2]$; 
observe that $x_{\delta}>c$ by the minimality of $\bar\xi$, and $(x_\delta,t_\delta)\to (c,t_0)$ as $\delta\to 0^+$. 
By \eqref{ad0-}, for such values of $\delta$ there holds $\psi_t(x_\delta,t_\delta)+H(\psi_x(x_\delta,t_\delta)-\delta)\le 0$. Letting $\delta\to 0^+$ we obtain that $\varphi_t(c,t_0)+H(\bar\xi)\le 0$. Since $\bar\xi\le\varphi_x(c^+,t_0)$, we obtain \eqref{adter}.
\smallskip

Now let $\bar\xi=-\infty$. Then there exists a sequence $\xi_n\to -\infty$ such that \eqref{x2} holds for $\xi=\xi_n$, thus for all $\xi\geq \xi_n$. 
By the arbitariness of $n$, \eqref{x2} holds for all $\xi\le \varphi_x(c,t_0)$. 
Hence for any $\xi\le \varphi_x(c,t_0)$ there exists $\psi\in C^1(\hat{Q}_1)$ 
such that $\psi_x(c,t_0)=\xi$, $\psi_t(c,t_0)=\varphi_t(c,t_0)$, and $u^*-\psi$ has a  strict maximum at $(c,t_0)$. 

Since, by Lemma \ref{propis}, 
$
\displaystyle{u^*(c,t_0)= \!\!\!\limsup_{\overline{Q}_1\ni (y,\tau)\to (c,t_0)} u^*(y,\tau)={\rm ess}\!\!\!\!\!\limsup_{Q_1\ni (y,\tau)\to (c,t_0)}u^*(y,\tau),}
$
there exists $\{(y_n,\tau_n)\}\subseteq \overline{Q}_1$ such that 
\begin{equation}\label{yn prop}
y_n>c\,,\quad\ (y_n,\tau_n)\to (c,t_0)\,,\quad\ u^*(y_n,\tau_n)\to u^*(c,t_0)\,.
\end{equation} 
In particular, for every $\ep>0$ there exists $n_{\ep}\in\mathbb{N}$ such that
\begin{equation}\label{st0}
u^*(c,t_0)-\ep < u^*(y_n,\tau_n)<u^*(c,t_0)+\ep\quad\mbox{for all}\ \,n>n_{\ep}\,.
\end{equation}
Let $\ep>0$ be fixed, and for $n>n_{\ep}$ set
\begin{equation}\label{st00}
\psi_{\ep,n}(x,t)=\psi(x,t)-2\ep \left(1- \left[1-\frac{x-c}{y_n-c}\right]_+^2\right)\qquad ((x,t)\in\hat{Q}_1)\,.
\end{equation}
Then $\psi_{\ep,n}\in C^1(\hat{Q}_1)$, and 
\begin{equation}\label{st1}
\left(\psi_{\ep,n}\right)_t=\psi_t\,,\quad \left(\psi_{\ep,n}\right)_x = \psi_x-4\ep\,\frac{[y_n-x]_+}{(y_n-c)^2}  \quad\mbox{in}\ \,\hat{Q}_1\,,
\end{equation}
\begin{equation}\label{st2}
\psi_{\ep,n}(c,t)=\psi(c,t)\ \ \mbox{for all $t\in(t_1,t_2]$}\,,
\end{equation}
\begin{equation}\label{st3}
 \psi-2\ep \leq \psi_{\ep,n}\leq \psi \quad\mbox{in}\ \,\hat{Q}_1\,.
\end{equation}

Without loss of generality we may assume that $t_0<t_2$. 
We fix a positive constant $\sigma <\min\{ d-c,t_0-t_1,t_2-t_0\}$ and set $Q_2=(c,c+\s)\times(t_0-\s,t_0+\s)$.  
\smallskip

\noindent {\it Claim:} If $\ep>0$ is sufficiently small and $n>n_\ep$, then $u^*-\psi_{\ep,n}$ has a
maximum in $\hat{Q}_1$ which is attained at a point $(x_{n,\ep},t_{n,\ep})=(x_n,t_n)\in Q_2$. 
\smallskip

In fact, since $u^*-\psi$ has a  strict maximum at $(c,t_0)$ and \eqref{st3} holds, for all $\ep>0$ small enough and $n>n_\ep$ the function 
$u^*-\psi_{\ep,n}$ has a maximum in $\hat{Q}_1$, 
attained at a point $(x_n,t_n)\in[c,c+\s)\times(t_0-\s,t_0+\s)$.
It remains to prove that $x_n>c$. Suppose that $x_n=c$. Then $t_n=t_0$, since $u^*-\psi$ has a 
strict maximum at $(c,t_0)$ and $\psi_{\ep,n}(c,t)=\psi(c,t)$ for every $t\in(t_0-\s,t_0+\s)$ (see \eqref{st2}). On the other hand, by \eqref{st0} and the equality $\psi_{\ep,n}(y_n,\tau_n)=\psi(y_n,\tau_n)-2\ep$ (see \eqref{st00}), we have that
\begin{equation*}
u^*(y_n,\tau_n)-\psi_{\ep,n}(y_n,\tau_n) \,=\, u^*(y_n,\tau_n) - \left [\psi(y_n,\tau_n)-2\ep \right]> u^*(c,t_0) -\psi(y_n,\tau_n) +\ep\,.
 \end{equation*}
Since $\psi(y_n,\tau_n)\to \psi(c,t_0)$ (see \eqref{yn prop}), for all sufficiently large $n$ we also get
\begin{equation*}
u^*(y_n,\tau_n)-\psi_{\ep,n}(y_n,\tau_n)> u^*(c,t_0) -\psi(c,t_0)= u^*(c,t_0) -\psi_{\ep,n}(c,t_0).
\end{equation*}
Since $y_n>c$ by \eqref{yn prop}, we have found a contradiction and proved the Claim.

\smallskip

Since $u$ is a viscosity subsolution of  $u_t+H(u_{x})=0$  in $Q_1$, it follows from the Claim, 
\eqref{ad0-} and the first equality in \eqref{st1}, that   
\begin{equation}\label{st4}
(\psi_{\ep,n})_t(x_n,t_n)+H((\psi_{\ep,n})_x(x_n,t_n))=\psi_t(x_n,t_n)+H((\psi_{\ep,n})_x(x_n,t_n))\leq 0\,.
 \end{equation}
Since $(\psi_{\ep,n})_x(x_n,t_n)\leq \psi_x(x_n,t_n)$ (see \eqref{st1}), it follows from \eqref{st4} that
\begin{equation}\label{st5}
\psi_t(x_n,t_n)+\inf\{H(s)|\,s\leq \psi_x(x_n,t_n)\}\leq \psi_t(x_n,t_n)+H((\psi_{\ep,n})_x(x_n,t_n))\leq 0\,.
 \end{equation}
On the other hand, letting first $n\to \infty$ and then $\ep\to 0^+$ in \eqref{st3} we find that $(x_n,t_n)=(x_{n,\ep},t_{n,\ep})\to(c,t_0)$ (recall that $c\leq x_{n,\ep}\leq c+\sigma$ and $\sigma>0$ is arbitrarily fixed). Since $\psi_x(c,t_0)=\xi$, we obtain from \eqref{st5}  that
$$
\psi_t(c,t_0)+\inf\{H(s)|\,s\leq \xi\}\leq 0 \quad\text{for all $\xi\in (-\infty,\varphi_x(c,t_0)]$.}
$$
Hence inequality \eqref{adter} also holds if $\bar\xi=-\infty$. 
\end{proof}

\begin{prop}\label{dominflubis} Let $H\in {\rm Lip}(\R)$. 
Let  the trapezoid
$A$ be defined by \eqref{trap bis} and let $\varphi\in C^1( \hat{A})$, where 
$\hat{A}\,:=\,\overline A\cap (\R\times (0,T])$.

\noindent $(i)$ Let $u$ be a viscosity subsolution of equation \eqref{hje} in $A$.  
Let $\big(d-\|H'\|_{\infty}t_0,t_0\big)$, with $t_0\in(0, T]$, be a local maximum point of $u^*-\varphi$ in $\hat{A}$. Then
\begin{equation}\label{dinf3}
\varphi_t\big(d-\|H'\|_{\infty}t_0,t_0\big)+H\big(\varphi_x\big((d-\|H'\|_{\infty}t_0)^-,t_0\big)\big)\le 0\,.
\end{equation} 

\noindent $(ii)$ Let $v$ be a viscosity supersolution of equation \eqref{hje} in $A$. Let $\big(d-\|H'\|_{\infty}t_0,t_0\big)$, $t_0\in (0,T]$, be a local minimum point of 
$v_*-\varphi$ in $\hat{A}$. Then
\begin{equation}\label{dinf4}
\varphi_t\big(d-\|H'\|_{\infty}t_0,t_0\big)+H\big(\varphi_x\big((d-\|H'\|_{\infty}t_0)^-,t_0\big)\big)\ge 0\,.
\end{equation} 
Similar results hold for  the trapezoids $B$ and $C$.
\end{prop}

\begin{proof} We only prove claim $(i)$ for $A$.  The remaining proofs are similar. Set 
$$               
\text{$ w(y,t)=u(y-\|H'\|_{\infty}t,t)$ \quad for $(y,t)\in \mathcal{A}=(a+\|H'\|_{\infty}T,d)\times (0,T)$} 
$$
(observe that under the inverse transformation $(y,t)\mapsto(x,t)=(y-\|H'\|_{\infty}t,t)$ the set $\mathcal{A}$ is mapped onto $\{(x,t) \,|\,x\in
(a+\|H'\|_{\infty}(T-t), d-\|H'\|_{\infty}t), \,t\in(0,T)\}$, which is a proper subset  of $A$).

We set $\tilde{H}(p)=H(p)+\|H'\|_\infty p$ for  $p\in\R$. Then $\tilde H$
is Lipschitz continuous on $\R$ and nondecreasing. 
We claim that $ w$ is a viscosity subsolution in  $\mathcal{A}$ of the equation 
\begin{equation}\label{eqV}
 w_t+\tilde{H}( w_y)=0\quad\text{in }\mathcal{A}\,.
\end{equation}
In fact, fix any $\psi=\psi(y,t)\in C^1\left(\hat{\mathcal{A}}\right)$, where $\hat{\mathcal{A}}=\overline {\mathcal A}\cap(\R\times (0,T])$,  and let 
$w^*-\psi$ have a local maximum at $(\bar y,\bar t)\in \mathcal{A}$; here by definition (see \eqref{use})
$$
 w^*(y,t)=\text{\rm ess}\!\!\!\!\!\!\limsup_{\mathcal{A}\ni (\eta,\t)\to (y,t)}\!\!  w(\eta,\t) \quad\text{for any }(y,t)\in \overline{\mathcal{A}}.
$$ 
It is easily seen that $ w^*(y,t)=u^*(y-\|H'\|_{\infty}t,t)$ for every $(y,t)\in \overline{\mathcal{A}}$ such that $a+\|H'\|_{\infty}T<y\leq d$. 
Therefore $u^*-\varphi$, where $\varphi=\varphi(x,t)=\psi(x+\|H'\|_{\infty}t,t)$,  has a local maximum at $(\bar{y}-\|H'\|_{\infty}\bar{t},\bar{t})\in A$. 
Since $u$ is a viscosity subsolution of equation \eqref{hje} in $A$, by \eqref{ad0-}, we obtain the claim:
$$
\begin{aligned}
&\varphi_t(\bar{y}-\|H'\|_{\infty}\bar{t},\bar{t})+H(\varphi_x(\bar{y}-\|H'\|_{\infty}\bar{t},\bar{t}))\,=\\
&\quad =\psi_t(\bar{y},\bar{t})+\|H'\|_\infty\psi_y(\bar{y},\bar{t})+H(\psi_y(\bar{y},\bar{t}))
=\psi_t(\bar{y},\bar{t})+\tilde{H}(\psi_y(\bar{y},\bar{t}))\le 0 \,.
\end{aligned}
$$

Now let $\varphi\in C^1(\hat{A})$, and let $(d-\|H'\|_{\infty}t_0,t_0)$ $(t_0\in(0, T])$ be a local maximum point of $u^*-\varphi$ in $\hat{A}$. Then $(d,t_0)$ is 
a local maximum point of $ w^*-\psi$ in $\hat{\mathcal{A}}$, where $\psi\in C^1\left(\hat{\mathcal{A}}\right)$ is defined by $\psi=\psi(y,t)=\varphi(y-\|H'\|_{\infty}t,t)$. 
Since $w$ is a viscosity subsolution in  $\mathcal{A}$ of \eqref{eqV} and $\tilde{H}$ is nondecreasing, by \eqref{adterbis} we have that
$$
\begin{aligned}
0&\geq\, \psi_t(d,t_0)+\inf\{\tilde{H}(\xi)\,|\,\xi\ge \psi_y(d^-,t_0)\} =\psi_t(d,t_0)+\tilde{H}(\psi_y(d^-,t_0))\,=\\
&=\psi_t(d,t_0)+\|H'\|_{\infty}\psi_y(d^-,t_0)+H(\psi_y(d^-,t_0))\,=\\
&=\varphi_t\big(d-\|H'\|_{\infty}t_0,t_0\big)+H\big(\varphi_x\big((d-\|H'\|_{\infty}t_0)^-,t_0\big)\big)\,. 
\end{aligned}
$$
\end{proof}

Now we can prove Theorem \ref{comp} for trapezoidal domains. We give the proof for problem $(T_{a,\pm})$ in the domain $A$.
The proof is similar for problem $(T^{a,\pm})$ in $B$ and easier for problem $(T_C)$.
\medskip

\noindent {\em Proof of Theorem \ref{comp} for problem $(T_{a,\pm})$.} 
Arguing by contradiction we suppose that  
\begin{equation}\label{sigma}
\max_{\overline  A}\, [u^*-v_*]_+ = \max_{[a,d]}\, \left[u^*(\cdot,0)-v_*(\cdot,0)\right]_++\sigma
\end{equation} 
for some $\sigma>0$. Consider the function $ F:  \overline A\times \overline A\mapsto\R$ defined by
\begin{equation}\label{Fipm}
F(x,t,y,s)\,:=\, u^*(x,t)- v_*(y,s)-\lambda (t+s)-\frac{|\,x-y + \a\ep\,|^p+ |\,t-s\,|^p}{\ep^p}\,, 
\end{equation}
where $\l\in\left(0, \frac{\s}{4T}\right)$ is fixed, $p\in(1,2]$ and $\ep\in(0,\min\{d-a,1\})$ with
\begin{equation}\label{ep0}
\a=\pm\tfrac 12\sqrt\ep\qquad \text{(here $\pm$ refers to problem } ({\rm T}_{a,\pm}). 
\end{equation} 
This implies that $|\,\a|\ep\le\frac{d-a}{2}$, and so $\frac{a+d}{2}+ \a\ep\in[a,d]$.

Since $F$ is upper semicontinuous,  it attains the maximum in $ \overline A\times \overline A$ at some point $( \bar x,\bar t, \bar y, \bar s)$. 
Observe that $F$ and $( \bar x,\bar t, \bar y, \bar s)$ depend on $\ep$ and $p$, but for notational simplicity we suppress this dependence. 

Since $\frac{a+d}{2}+ \a\ep\in[a,d]$, there holds 
$$
   F\left(\tfrac12 (a+d),0,\tfrac12 (a+d)+ \a\ep,0\right)\le  F\!\left( \bar x,\bar t, \bar y, \bar s\right) \,,
$$
which implies that 
$$
\frac{|\, \bar x- \bar y + \a\ep\,|^p+|\,\bar t- \bar s\,|^p}{\ep^p}\le 4M, \quad \text{where $M=\max\left\{\| u^*\|_{L^\infty( A)},\| v_*\|_{L^\infty( A)}\right\}$.}
$$
Hence there holds $|\,\bar t- \bar s \,|\,\le \,(4M)^{\frac1p}\,\ep$ and
\begin{equation}\label{eps1}
|\, \bar x- \bar y+ \a\ep\,|\,\le \,(4M)^{\frac1p}\,\ep \quad \Rightarrow \quad | \bar x- \bar y|\,\le \, 
\left[(4M)^{\frac1p}+\tfrac12\right]\ep. 
\end{equation} 
Observe that both estimates can be made independent of $p\in (1, 2]$. 

Now consider the function $g:  \overline A\mapsto\R$,
\begin{equation}\label{gipm}
g=  g(x,t)\,:=\, u^*(x,t)- v_*(x,t)-2\lambda t =  F(x,t,x,t) +|\,\a|^p\,.
\end{equation}
Observe that $g$ is upper semicontinuous, thus its maximum in $\overline A$ exists.

Set $ A_\t=\{(x,t)\in \overline A\,|\, t\in[0,\t]\}$ for any $\t\in(0,T]$. 
Below we prove the following claims. 
\smallskip

\noindent {\em Claim 1:} There exists $\tau\in (0,T)$ such that
\begin{equation}\label{cl1}
\max_{ \overline A} g= \max_{ \overline  A\setminus A_\t} g \quad\mbox{and}\quad \max_{\overline A} g>\max_{A_\t} g\,.
\end{equation}
{\em Claim 2:} There exists $\ep_1\in(0,d-a)$ which does not depend on $p\in (1,2]$ 
such that for all $\ep\in(0,\ep_1)$ and $p\in(1,2]$ 
\begin{equation}\label{cl2}
\max_{ \overline A\times \overline  A} F= \max_{\left( \overline A\setminus A_\t\right)^2} F, \qquad \text{with $\tau\in (0,T)$ given by Claim 1.}
\end{equation}

To prove Claim 1 set $G(t)\,:=\,\max\limits_{A_t} g$\, $(t\in(0,T])$.
By \eqref{sigma} and since  $\l< \frac{\s}{4T}$, 
\begin{eqnarray}\label{cl11}
&& G(T) \ge\max_{ \overline A}\,(u^*-v_*)-2\lambda T=\max_{ \overline A}\,[u^*-v_*]_+-2\lambda T\,=\\
&&=\,  \max_{[a,d]}\, \left[u^*(\cdot,0)-v_*(\cdot,0)\right]_++\sigma-2\lambda T\,\ge
\max_{[a,d]}\, \left[u^*(\cdot,0)-v_*(\cdot,0)\right]_++ \tfrac12 \s\,.\nonumber
\end{eqnarray}
Since $G$ is nondecreasing, there exists $L_0\,:=\,\lim\limits_{t\to 0^+}G(t)$, and Claim 1 follows if we prove that
\begin{equation}\label{cl12}
L_0\le \max_{[a,d]}\, \left[u^*(\cdot,0)-v_*(\cdot,0)\right]_+\,.
\end{equation}
In fact, by \eqref{cl11}-\eqref{cl12} there exists $\tau\in (0,T)$ such that
\begin{equation}\label{cl13}
G(\t)=\max_{ A_\t}g\le \max_{[a,d]}\, \left[u^*(\cdot,0)-v_*(\cdot,0)\right]_++\tfrac14\s
<G(T)=\max_{\overline A}g\,.
\end{equation}

To prove \eqref{cl12}, let $\{t_n\}$ be a decreasing sequence such that $ t_n\to 0^+$, and let $(x_n,\tau_n)\in A_{t_n}$ be a  maximum point  - namely, $G(t_n)=g(x_n,\t_n)$. Clearly, there exists a converging subsequence (not relabelled) of $\{(x_n,\tau_n)\}\subseteq A_{t_1}$  and a point $(\overline x,0)$,  $\overline x\in [a,d]$, such that $(x_n,\tau_n)\to(\overline x,0)$ as $n\to\infty$ (observe that $\lim_{n\to \infty} \tau_n=\lim_{n\to \infty} t_n=0$). Then, by the upper semicontinuity of $g$,
$$
\begin{aligned}
L_0&=\lim_{n\to \infty}G(t_n) =\lim_{n\to \infty}\max\limits_{A_{t_n}}g\le g(\overline x,0^+)\,=
u^*(\overline x,0^+)- v_*(\overline x,0^+)\,\leq\\
&\le  \max_{[a,d]}\, \left[u^*(\cdot,0)-v_*(\cdot,0)\right]_+\,.
\end{aligned}
$$
This proves \eqref{cl12} and Claim 1 follows.

\smallskip

To prove Claim 2, we preliminarily observe that, by \eqref{Fipm} and \eqref{gipm}, 
for every maximum point $(\bar{x},\bar{t}, \bar{y}, \bar{s})$ of $F$  there holds
\begin{equation*}
\max_{\overline A}\, g -|\,\a|^p \,\le\, F(\bar{x},\bar{t}, \bar{y}, \bar{s})\,\le\, u^*(\bar{x},\bar{t})-v_*(\bar{y},\bar{s})-\lambda(\bar{t}+\bar{s})\,.
\end{equation*}
In view of \eqref{ep0}, this implies that 
\begin{equation}\label{dest2}
\max_{\overline  A} g -\sqrt\ep \,\le\, u^*(\bar{x},\bar{t})-v_*(\bar{y},\bar{s})-\lambda(\bar{t}+\bar{s}).
\end{equation}

Now we argue by contradiction. 
Were Claim 2 false, there would exist a sequence $\{\ep_n\}\subset (0,\ep_0)$ such that $ \ep_n\to 0^+$, a sequence $\{p_n\}\subset (1,2]$
and a sequence of maximum points of $F$
$$
\{(\bar{x}_n,\bar{t}_n, \bar{y}_n, \bar{s}_n)\}\equiv\{(\bar{x}(\ep_n,p_n),\bar t(\ep_n,p_n), \bar y(\ep_n,p_n), \bar s(\ep_n,p_n))\}\subset A_\t\,.
$$ 
By the boundedness of $\{(\bar{x}_n,\bar{t}_n, \bar{y}_n, \bar{s}_n)\}$ and 
\eqref{eps1}, there would exist a converging subsequence (not relabelled) of $\{(\bar{x}_n,\bar{t}_n, \bar{y}_n, \bar{s}_n)\}$  
and a point $(\tilde{x},\tilde{t},\tilde{x}, \tilde{t})\in \left( A_\t\right)^2$, such that 
$(\bar{x}_n,\bar{t}_n, \bar{y}_n, \bar{s}_n)\to(\tilde{x},\tilde{t},\tilde{x}, \tilde{t})$ as $n\to\infty$. 
Rewriting \eqref{dest2} with 
$\ep=\ep_n$, $(\bar{x},\bar{t}, \bar{y}, \bar{s})=(\bar{x}_n,\bar{t}_n, \bar{y}_n, \bar{s}_n)$ and letting $n\to \infty$, it follows 
from the upper semicontinuity of the function at the right-hand side of \eqref{dest2} that 
$\max\limits_{\overline A} g\le g(\tilde{x},\tilde{t})$,
which contradicts Claim 1 since $(\tilde{x},\tilde{t})\in  A_\t$. 
Hence we have also proved Claim 2.

\smallskip

Now we complete the proof. Henceforth we assume that $\ep\in(0,\ep_1)$, so that we can use Claim 2. Then the function 
 \begin{eqnarray*}
(x,t)\mapsto F(x,t, \bar{y}, \bar{s})\!\!\!\!&=& \!\!\!\! u^*(x,t)- v_*( \bar{y}, \bar{s})-\lambda (t+ \bar{s}) -\frac{|\,x- \bar{y} + \a\ep\,|^p+|\,t- \bar{s}\,|^p}{\ep^p}\,=: \\
&=:&\!\!\!\! u^*(x,t)-\phi(x,t)
 \end{eqnarray*}
has a maximum at some point $(\bar{x},\bar{t})\in  \overline A\setminus A_\t$ with $\tau\in (0,T)$ given by Claim 1. Similarly, the function
\begin{eqnarray*}
(y,s)\mapsto
- F( \bar{x},\bar{t},y,s)\!\!\!\!&=& \!\!\!\!   v_*(y,s)- u^*(\bar{x},\bar{t})+\lambda ( \bar{t}+s)+\frac{|\,\bar{x}-y + \a\ep\,|^p+| \,\bar{t}-s\,|^p}{\ep^p}\,=: \\
&=& \!\!\!\!  v_*(y,s)-\chi(y,s)
\end{eqnarray*}
has a minimum at some point $( \bar{y}, \bar{s})\in  \overline A\setminus A_\t$ with $\tau\in (0,T)$ as above. 
Since $u$ is a viscosity subsolution and $v$ a viscosity supersolution of $({\rm T}_{a,\pm})$, 
by definition $\phi$ and $\chi$ must satisfy suitable differential inequalities at $(\bar{x},\bar{t})$ and $( \bar{y}, \bar{s})$ 
(see Subsection \ref{def} and Proposition \ref{p2}). We show below that these inequalities always lead to a contradiction, whence the result follows. 
Before proceeding observe that
$$
\phi_t(x,t)=\l +\frac{p\,|\,t-\bar{s}\,|^{p-1}\,{\rm sgn}(t-\bar{s})}{\ep^p}\,,\quad
\phi_x(x,t)= \frac{p\,|\,x- \bar{y} + \a\ep\,|^{p-1}\,{\rm sgn}(x- \bar{y} + \a\ep)}{\ep^p}\,,
$$
$$
\chi_s(y,s)=-\l +\frac{p\,|\,\bar{t}-s\,|^{p-1}\,{\rm sgn}(\bar{t}-s)}{\ep^p}\,,\quad
\chi_y(y,s)= \frac{p\,|\,\bar{x}- y + \a\ep\,|^{p-1}\,{\rm sgn}(\bar{x}- y + \a\ep)}{\ep^p}\,.
$$

\medskip

\noindent {\em Problem }$({\rm T}_{a,+}).$  Then $\a=\tfrac12\sqrt\ep$ (see \eqref{ep0}) and  
we distinguish two cases: 

\smallskip

$(1)\; \bar{x}=\bar{x}_{\ep,p}\geq a,\,\,\bar{y}=\bar{y}_{\ep,p}>a$ for some $\ep\in (0,\ep_1)$ and $p\in (1,2 ]$: 
Since $u$ is a viscosity subsolution and $( \bar{x}, \bar{t})$ is a maximum point of $u^*-\phi$, it follows from \eqref{ad1-} and \eqref{ad0-} that  
\begin{eqnarray}\label{eq117}
&&\phi_t(\bar{x},\bar{t})+H(\phi_x(\bar{x},\bar{t}))\\ 
&&\quad =\l +\frac{p\,|\,\bar{t}-\bar{s}\,|^{p-1}\,{\rm sgn}(\bar{t}-\bar{s})}{\ep^p}
+H\left(\frac{p\,|\,\bar{x}- \bar{y} + \a\ep\,|^{p-1}\,{\rm sgn}(\bar{x}- \bar{y} + \a\ep)}{\ep^p}\right)\le 0\,. \nonumber
\end{eqnarray}
Similarly, since $v$ is a viscosity supersolution and $( \bar{y}, \bar{s})$ is a minimum point of $v_*-\chi$, it follows from \eqref{ad0+} (if $\overline y<d-\|H'\|_{\infty}\overline s$) and  (\ref{dinf4}) (if $\overline y=d-\|H'\|_{\infty}\overline s$) that
\begin{eqnarray}\label{eq118}
&&\chi_s(\bar{y},\bar{s})+H(\chi_y(\bar{y},\bar{s}))\\ 
&&\quad =-\l +\frac{p\,|\,\bar{t}-\bar{s}\,|^{p-1}\,{\rm sgn}(\bar{t}-\bar{s})}{\ep^p}
+H\left(\frac{p\,|\,\bar{x}- \bar{y} + \a\ep\,|^{p-1}\,{\rm sgn}(\bar{x}- \bar{y} + \a\ep)}{\ep^p}\right)\ge0\,. \nonumber
\end{eqnarray}
Subtracting \eqref{eq118} from \eqref{eq117} we find that $2\lambda\le 0$, which is a contradiction. 

\smallskip

$(2)\;\bar{x}=\bar{x}_p\geq a,\,\bar{y}=\bar{y}_p=a$ for all $p\in (1, 2 ]$ and $\ep\in (0,\ep_1)$:  we fix $\ep=\tilde{\ep}\in(0,\ep_1)$ so small that
 \begin{equation}\label{hyp+}
\sup\limits_{\xi\ge 1/\tilde{\ep}  }H(\xi)\,\le\,\limsup\limits_{\xi\to\infty}H(\xi)+ \tfrac12 \lambda \,,\quad \;
H\left(1/\tilde{\ep}\right)\,\ge\, \limsup\limits_{\xi\to\infty}H(\xi)-\tfrac12 \lambda \,.
\end{equation}
Since we have chosen $\ep=\tilde{\ep}$, $\bar{x}=\bar{x}_p$ and $\bar{y}=\bar{y}_p$ only depend on $p$. Now \eqref{eq117} reads
\begin{equation}\label{eq1177}
\phi_t(\bar{x},\bar{t})+H(\phi_x(\bar{x},\bar{t}))\,=
\l +\frac{p\,|\,\bar{t}-\bar{s}\,|^{p-1}\,{\rm sgn}(\bar{t}-\bar{s})}{\tilde{\ep}^p}
+H\left(\frac{p\,|\,\bar{x}- a + \a\tilde{\ep}\,|^{p-1}}{\tilde{\ep}^p}\right)\le 0
\end{equation}
(since $\bar{x}- a + \a\tilde{\ep}>0$, $\,{\rm sgn}(\bar{x}- a + \a\tilde{\ep})=1$). On the other hand, by \eqref{adqua} there holds
 \begin{eqnarray}\label{eq113}
&&\chi_s(a,\bar{s})+\sup\{H(\xi)\,|\,\xi\ge \chi_y(a^+,\bar{s})\}\\
&&\quad =\,-\l +\frac{p\,|\,\bar{t}-\bar{s}\,|^{p-1}\,{\rm sgn}(\bar{t}-\bar{s})}{\tilde{\ep}^p}+\sup\left\{H(\xi)\,\Big|\,\xi\ge \frac{p\,|\,\bar{x}- a + \a\tilde{\ep}\,|^{p-1}}{\tilde{\ep}^p} \right\}
\ge 0. \nonumber
\end{eqnarray}  
Subtracting \eqref{eq113} from \eqref{eq1177} we find that 
 \begin{equation*}
2\l+H\left(\frac{p\,|\,\bar{x}- a + \a\tilde{\ep}\,|^{p-1}}{\tilde{\ep}^p}\right)\,\le \sup\left\{H(\xi)\,\Big|\,\xi\ge \frac{p\,|\,\bar{x}- a + \a\tilde{\ep}\,|^{p-1}}{\tilde{\ep}^p} \right\}
\end{equation*} 
whence, as $p\to1^+$,  
 \begin{equation}\label{zz98}
2\l+H\left(1/\tilde{\ep}\right)\,\le \sup\left\{H(\xi)\,\Big|\,\xi\ge 1/\tilde{\ep} \right\}.
\end{equation}  
By \eqref{hyp+} this implies that  $\l\le0$, and again we have found a contradiction.

\medskip

\noindent {\em Problem }$({\rm T}_{a,-})$.  Then $\a=-\tfrac12\sqrt\ep$ (see \eqref{ep0}) and  
again we distinguish two cases: 

\smallskip

$(1)\; \bar{x}=\bar{x}_{\ep,p}> a,\,\,\bar{y}=\bar{y}_{\ep,p}\geq a$ for some $\ep\in (0,\ep_1)$ and $p\in (1,2 ]$: in this case 
\eqref{eq117}-\eqref{eq118} follow from \eqref{ad0-}, \eqref{ad0+} and \eqref{ad2+}, whence again $2\lambda\le 0$. 

\smallskip

$(2)\;\bar{x}=\bar{x}_p= a,\,\bar{y}=\bar{y}_p\geq a$ for all $p\in (1, 2 ]$ and $\ep\in (0,\ep_1)$:  we fix $\ep=\tilde{\ep}\in(0,\ep_1)$ so small that
\begin{equation}\label{hyp-}
\inf\limits_{\xi\le -1/\tilde{\ep}  }H(\xi)\,\ge\,\liminf\limits_{\xi\to-\infty}H(\xi)- \tfrac12 \lambda\,,\quad \;
H\left(-1/\tilde{\ep}\right)\,\le\, \liminf\limits_{\xi\to-\infty}H(\xi)+\tfrac12\lambda
\end{equation}
(so $\bar{x}=\bar{x}_p$ and $\bar{y}=\bar{y}_p$ only depend on $p$). By \eqref{adter} there holds
\begin{eqnarray}\label{eq44}
&&\phi_t(a,\bar{t})+\inf\{H(\xi)\,|\,\xi\le \phi_x(a^+,\bar{t})\}\,=\\
&&\quad =\,\l +\frac{p\,|\,\bar{t}-\bar{s}\,|^{p-1}\,{\rm sgn}(\bar{t}-\bar{s})}{\tilde{\ep}^p}+\inf\left\{H(\xi)\,\Big|\,\xi\le -\frac{p\,|\,a -\bar{y}+ \a\tilde{\ep}\,|^{p-1}}{\tilde{\ep}^p} \right\} \le 0 \nonumber
\end{eqnarray}  
(since $a-\bar{y} + \a\tilde{\ep}<0$, $\,{\rm sgn}(a-\bar{y} + \a\tilde{\ep})=-1$). On the other hand,  since $\bar{y}>a$,
\begin{eqnarray}\label{eq445566}
&&\chi_s(\bar{y},\bar{s})+H(\chi_y(\bar{y},\bar{s}))\,=\\ 
&&\quad =\,-\l +\frac{p\,|\,\bar{t}-\bar{s}\,|^{p-1}\,{\rm sgn}(\bar{t}-\bar{s})}{\tilde{\ep}^p}
+H\left(-\frac{p\,|\,a- \bar{y} + \a\tilde{\ep}\,|^{p-1}}{\tilde{\ep}^p}\right)\ge0\,. \nonumber
\end{eqnarray}
Subtracting \eqref{eq445566} from \eqref{eq44} gives
 \begin{equation*}
2\l+\inf\left\{H(\xi)\,\Big|\,\xi\le -\frac{p\,|\,a -\bar{y}+ \a\tilde{\ep}\,|^{p-1}}{\tilde{\ep}^p} \right\}\,\le 
H\left(-\frac{p\,|\,a- \bar{y} + \a\tilde{\ep}\,|^{p-1}}{\tilde{\ep}^p}\right)\,.
\end{equation*} 
Letting $p\to1^+$ this implies that 
 \begin{equation}\label{zz898}
2\l+\inf\left\{H(\xi)\,\Big|\,\xi\le -1/\tilde{\ep} \right\} \,\le H\left(-1/\tilde{\ep}\right) \,,
\end{equation}  
whence by \eqref{hyp-} we get $\l\le0$, again a contradiction.
\hfill$\square$

\smallskip

It remains to prove Theorem \ref{comp} for problem $(N_{a,\pm}^{b,\pm})$. 
\smallskip

\noindent {\em Proof of Theorem \ref{comp} for problem $(N_{a,\pm}^{b,\pm})$.} 
Let $Q=(a,b)\times(0,T)$ with $-\infty<a<b<\infty$, and $Q_\t=\{(x,t)\in Q\,|\, t\in[0,\t]\}$ for any $\t\in(0,T]$. 
Let $u^*,v_*$ be defined in $\overline{Q}$ by  \eqref{use}-\eqref{lse}.
As before we set
$$
A_\t\,:=\,\{(x,t)\,|\,x\!\in\![a,b\!-\!\|H'\|_{\infty}t],\,t\!\in\! [0,\tau ]\}, \quad
B_\t\,:=\,\{(x,t)\,|\,x\!\in\![a\!+\!\|H'\|_{\infty}t,b],\,t\!\in\! [0,\tau ]\}
$$
for all $\tau\in (0,\tau_1]$, with $\t_1$ defined by \eqref{det1}. 
By Remark \ref{remres}, the restrictions $u_{1,A}=u\lefthalfcup A_{\tau_1}$, $v_{1,A}=v\lefthalfcup A_{\tau_1}$ 
are viscosity sub- and supersolutions of problem $({\rm T}_{a,\pm})$ in $A_{\t_1}$
(similarly for $u_{1,B}=u\lefthalfcup B_{\tau_1}$, $v_{1,B}=v\lefthalfcup B_{\tau_1}$). 
Since we have already proved Theorem \ref{comp} for the trapezoidal domains $A$ and $B$, we have that 
\begin{equation*}
\max_{A_{\t_1}} [\!(u_{1,A})^*\!-(v_{1,A})_*]_+\! \le  \! \max_{[a,b]}\left[\!(u_{1,A})^*\!(\cdot,0)\!-\!(v_{1,A})_*(\cdot,0)\right]_+
\!\leq\! \max_{[a,b]} \left[u^*\!(\cdot,0)\!-\!v_*\!(\cdot,0)\right]_+
\end{equation*}
(notice that $(u_{1,A})^*(x,t)\leq u^*(x,t)$ and $(v_{1,A})_*(x,t)\geq v_*(x,t)$, since $A_{\t_1}\subseteq Q$), 
\begin{equation*}
\max_{B_{\t_1}}[\!(u_{1,B})^*\!-(u_{1,B})_*]_+\! \le  \! \max_{[a,b]} \left[\!(u_{1,B})^*(\cdot,0)\!-\!(v_{1,B})_*(\cdot,0)\right]_+
\!\le\!   \max_{[a,b]} \left[u^*(\cdot,0)\!-\!v_*(\cdot,0)\!\right]_+\!,
\end{equation*}
 whence, by Remark \ref{restri},
$\displaystyle{\max_{[a,b]\times [0,\tau_1-\d]}\, [u^*-v_*]_+  \le  \, \max_{[a,b]}\, \left[u^*(\cdot,0)-v_*(\cdot,0)\right]_+}\,\,$
for all $\delta \in (0,\tau_1)$. By the arbitrariness of $\d$ this means that 
\begin{equation}\label{step1}
\sup_{[a,b]\times [0,\tau_1)}\, [u^*-v_*]_+ \le  \, \max_{[a,b]}\, \left[u^*(\cdot,0)-v_*(\cdot,0)\right]_+\,.
\end{equation}

Let $\d \in (0,\tau_1)$ be arbitrary and fixed. 
Arguing as before in the rectangle $\overline{Q}_{\t_2-\d}\setminus \overline{Q}_{\t_1-\d}$, where
$\t_2=\min\left\{(b-a)/\|H'\|_\infty,\, T\right\}$,
we obtain 
$$
[u^*(x,t)-v_*(x,t)]_+\leq \max_{[a,b]}\, \left[u^*(\cdot,(\tau_1-\d)^+)-v_*(\cdot,(\tau_1-\d)^+)\right]_+
$$
for all $(x,t)\in [a,b ]\times(\tau_1-\d,\tau_2-\d)$. Since 
$$
u^*(\cdot,(\t_1-\d)^+)-v_*(\cdot,(\t_1-\d)^+)\le u^*(\cdot,\t_1-\d)-v_*(\cdot,\t_1-\d)
$$ 
(see \eqref{lineq1}), from the above inequality and \eqref{step1} we obtain that
\begin{equation*}
\sup_{ [a,b ]\times[0,\tau_2-\d) }\, [u^*-v_*]_+\, \le  \, \max_{[a,b]}\, \left[u^*(\cdot,0)-v_*(\cdot,0)\right]_+\,,
\end{equation*}
whence, by the arbitrariness of $\d$,
\begin{equation*}
\sup_{ [a,b ]\times[0,\tau_2) }\, [u^*-v_*]_+\, \le  \, \max_{[a,b]}\, \left[u^*(\cdot,0)-v_*(\cdot,0)\right]_+\,.
\end{equation*}
It is now clear that in a finite number of steps the claim follows. 
\hfill$\square$


\section{Regularity}\label{Sec reg}
\setcounter{equation}{0}
In this section we prove Proposition \ref{prolip}. 
We begin with the proof of the first part, which concerns one-sided Lipschitz bounds in $t$ for 
sub- and supersolution of equation \eqref{hje}.
\medskip

\noindent {\em Proof of Proposition \ref{prolip}\,$(i)$.} 
To prove \eqref{sih}$_1$, i.e.\ the first inequality in \eqref{sih}, it is enough to show that  
\begin{equation}\label{Li0}
u^*(\bar x,t)\le u^*(\bar x,t_1^+)+K(t-t_1) 
\quad\text{for all }t\in (t_1,T],\,  \bar x\in\O=(a,b),
\end{equation}
Indeed, taking the $\limsup$ as $\bar x\to a^+$ (if $a\in \R$) or $\bar x\to b^-$ (if $b\in \R$) it follows from \eqref{tilque} that  \eqref{Li0} is also satisfied if $\bar x\in \partial\O$, and then \eqref{sih}$_1$ follows from \eqref{lineq1}.

We prove \eqref{Li0}. By Lemma \ref{propis} applied to the restriction of $u$ to $(a,b)\times (t_1,T)$, for any $\ep>0$ there exists  $\delta>0$ such that
\begin{equation}\label{Li1}
u^*(x,t_1^+) \le u^*(\bar{x},t_1^+) +\ep      \quad\text{for all }x\in (\bar{x}-\delta, \bar{x}+\delta)\,
\end{equation}
 and $(\bar{x}-\delta, \bar{x}+\delta)\subset \Omega$. Setting $D_\delta\,:=\,(\bar{x}-\delta, \bar{x}+\delta)\times(t_1,T)$ and $u_1\,:=\,u\lefthalfcup D_\delta$, it follows from Remark \ref{restri} that 
\begin{equation*}\label{eq se}
(u_1)^*(x,t)=
\begin{cases}
u^*(x,t)&\mbox{if}\ |x-\bar x|< \delta, \ t\in (t_1,T ]\,,\\ 
u^*(x,t_1^+)&\mbox{if}\ |x-\bar x|< \delta, \ t=t_1	\,,
\end{cases}
\end{equation*}
and from  Remark \ref{remres} that $u_1$ is a viscosity subsolution of equation \eqref{hje} in $D_\delta$. Observe that, by Definition \ref{defsubsupervi}$-(i)$, this implies that $u_1$ is also a viscosity subsolution of problem $(N_{\bar{x}-\delta, -}^{\bar{x}+\delta,+})$ in $D_\delta$.

Let 
\begin{equation}\label{Li2}
w(x,t)\,:=\,u^*(\bar{x},t_1^+) +\ep+ K(t-t_1)   \quad \text{if } (x,t)\in [\bar{x}-\delta, \bar{x}+\delta]\times[t_1,T].
\end{equation}
It is easy to prove that $w$ is a viscosity supersolution of problem $(N_{\bar{x}-\delta, -}^{\bar{x}+\delta,+})$ in 
$D_\delta$: if $\varphi \in C^1(\overline D_\delta)$ and $(x_0,t_0)\in \overline D_\delta$ (with $t_0>t_1$)  is a local minimum point of $w_*-\varphi =w-\varphi$ in $\overline D_\delta$, then $\varphi_t(x_0,t_0)+H(\varphi_x(x_0,t_0))\ge w_t(x_0,t_0)-K=0$.

Applying Theorem \ref{comp} in $D_\delta$, and observing that, by \eqref{Li1}-\eqref{Li2}, $(u_1)^*(\cdot,t_1)=u^*(\cdot,t_1^+)\le  w (\cdot,t_1)$ in $[\bar x-\delta, \bar x+\delta]$, we obtain that $u^*\le  w $ in $D_\delta$. In particular,  
$$
u^*(x,t)- u^*(\bar{x},t_1^+) \le \ep+K(t-t_1) \;\;\text{for any $|x-\bar{x}|<\delta$ and $t\in (t_1,T]$}
$$
and \eqref{Li0} follows from the arbitrariness of $\ep$.

The proof of \eqref{sih}$_2$ is similar: arguing as before one shows that
\begin{equation}\label{Li0b}
v_*(\bar x,t)\ge  v_*(\bar{x},t_1^+) +k(t-t_1) \;\;\text{for }t\in (t_1,T].
\end{equation}
\hfill$\square$

\smallskip

Formulas \eqref {wt 1}  and \eqref {wt 2} in Proposition \ref{prolip}$-(ii)$ quantify the fact that initial jump discontinuities cannot disappear instantaneously. They are a special case of the following result, with $x_0=c$, $t_0=0$ and $G=u_0$: 

\begin{lem}\label {ledi} Let $(H_1)$ hold, let $K,k$ be defined by \eqref{Def K,k}, let $u$ be a viscosity solution of equation \eqref{hje} and let $(x_0,t_0)\in \O\times [0,T)$. Let $G\in L^{\infty}_{{\rm loc}}(\overline{\Omega})$ be such that
\begin{equation}\label{ct01}
u^*(x,t_0^+)=G^*(x)\,,\quad\ u_*(x,t_0^+)=G_*(x)\,,
\end{equation}
\begin{equation}\label{hyp}
\text{there exist }G(x_0^{\pm})\,:=\,{\rm ess}\lim_{x\to x_0^{\pm}}G(x), \qquad 
G(x_0^+)\neq G(x_0^-),
\end{equation}
and set 
$$
\underline t:=\,\left\{\begin{array}{ll}
t_0+\frac{|G(x_0^+)- G(x_0^-)|}{K-k}>0&\mbox{if}\ \,K>k\,,\\
T&\mbox{{\rm otherwise}}\,.
\end{array}\right.
$$
Then for all $t\in (t_0,\min\{\underline t,T\})$
\begin{equation}\label{cl21}
G(x_0^+)> G(x_0^-)\ \Rightarrow \ 
\begin{cases}
u^*(x_0,t)=u^*(x_0^+,t)>u^*(x_0^-,t)\\
u_*(x_0^+,t)>u_*(x_0^-,t)=u_*(x_0,t),
\end{cases}
\end{equation}
\begin{equation}\label{cl22}
G(x_0^+)< G(x_0^-)\ \Rightarrow \ 
\begin{cases}
u^*(x_0^+,t)<u^*(x_0^-,t)=u^*(x_0,t)\\
u_*(x_0,t)=u_*(x_0^+,t)<u_*(x_0^-,t).
\end{cases}
\end{equation}
\end{lem}
\begin{proof} We only prove \eqref{cl21}, the proof of \eqref{cl22} is similar. 
By assumption, for any $\ep>0$ there exists $\delta\in(0,\bar{\d})$ such that
$$
|G(x)-G(x_0^-)|<\ep\ \ \mbox{for a.e.~$x\in (x_0-\delta,x_0)$}
$$
$$ 
|G(y)-G(x_0^+)|<\ep\ \ \mbox{for a.e.~$y\in (x_0,x_0+\delta)$}\,,
$$
whence, since $G(x_0^+)> G(x_0^-)$,
\begin{equation}\label{udx2}
G(x_0^-)-\ep \le G_*(x)\le G^*(x)\le G(x_0^-)+\ep\quad\text{for all $x\in(x_0-\delta,x_0)$}
\end{equation}
\begin{equation}\label{udx1}
G(x_0^+)-\ep \le G_*(y)\le G^*(y)\le G(x_0^+)+\ep\quad\text{for all $y\in(x_0,x_0+\delta)$.} 
\end{equation}
Then we get that for all $t\in(t_0,T ]$ and $x\in(x_0-\delta,x_0)$
$$
\begin{aligned}
&G(x_0^-)-\ep +k(t-t_0) \! \stackrel{\eqref{udx2}} \leq G_*(x)+k(t-t_0)\!\stackrel{\eqref{ct01}}=
u_*(x, t_0^+)+k(t-t_0)\!\stackrel{\eqref{Li0b}}\leq 
u_*(x,t)\le 
\\&\quad \le  
u^*(x,t)\!\stackrel{\eqref{Li0}} \leq \!  
u^*(x, t_0^+)\! +\! K(t\! -\! t_0)\!\stackrel{\eqref{ct01}}=
\! G^*(x)\! +\! K(t\! -\! t_0)\!\stackrel{\eqref{udx2}}  \leq \!  G(x_0^-)\! +\! \ep+K(t\! -\! t_0)\,. 
\end{aligned}
$$
Similarly, using \eqref{udx1} instead of \eqref{udx2}, for all $t\in (t_0,T ]$ and $y\in(x_0,x_0+\delta)$
$$
\begin{aligned}
&G(x_0^+)-\ep +k(t-t_0) \leq G_*(y)+k(t-t_0)
=u_*(y, t_0^+)+k(t-t_0)\leq \, u_*(y,t) 
\le \\&\quad 
\le u^*(y,t) \leq \, u^*(y, t_0^+)+K(t-t_0)=G^*(y)+K(t-t_0)\leq \, G(x_0^+)+\ep+K(t-t_0) \,.
\end{aligned}
$$
In particular we obtain from the above inequalities that for all $t\in (t_0,T ]$
\begin{subequations}\label{12}
\begin{equation}\label{1}
u_*(x,t)\le u^*(x,t) \le G(x_0^-)+\ep+K(t-t_0) \quad\text{for all $x\in(x_0-\delta,x_0)$, }
\end{equation}
\begin{equation}\label{2}
G(x_0^+)-\ep +k(t-t_0) \le u_*(y,t)\le u^*(y,t) \quad\text{for all $y\in(x_0, x_0+\delta)$. }
\end{equation}
\end{subequations}

Now set 
\begin{equation}\label{tbar}
\ep_0=\frac{G(x_0^+)-G(x_0^-)}{2}\,,\quad \underline t_\ep 
=\begin{cases}
t_0+\frac{2(\ep_0-\ep)}{K-k}  &\text{if $K>k$ \;$(\ep\in(0,\ep_0))$\,,} \smallskip\\ 
T&\text{otherwise.}
\end{cases}
\end{equation}
Observe that $\underline t_\ep\to \underline t$ as $\ep\to 0$. 
Then for all $\ep\in(0,\ep_0)$ and $t\in [t_0,\underline t_\ep)$ there holds
\begin{equation}\label{nuo} 
G(x_0^+)-\ep +k(t-t_0) \,>\, G(x_0^-)+\ep+K(t-t_0)\,,
\end{equation}
whence, by \eqref{12}, for all $\ep\in(0,\ep_0)$, $t\in (t_0,\underline t_\ep)$, $x\in (x_0-\delta,x_0)$ and $y\in (x_0,x_0+\delta)$
$$
u^*(x,t) \leq G(x_0^-)+\ep+K(t-t_0)<G(x_0^+)-\ep+k(t-t_0) \leq \, u^*(y,t)\,,
$$
$$
u_*(x,t)\leq  G(x_0^-)+\ep+K(t-t_0)<G(x_0^+)-\ep+k(t-t_0)\leq u_*(y,t) \,.
$$
Plainly this implies that for all $t\in (t_0,\underline t_\ep)$
\begin{eqnarray*}
&&u^*(x_0^-,t)={\rm ess}\!\!\!\!\!\!\!\!\limsup_{Q\ni(x,\tau)\to (x_0^-,t)}u^*(x,\tau)\,\leq G(x_0^-)+\ep+K(t-t_0)\,<\\
&&\qquad \qquad<\,G(x_0^+)-\ep+k(t-t_0)\,\leq\,{\rm ess}\!\!\!\!\!\!\!\!\limsup_{Q\ni(y,\tau)\to (x_0^+,t)}u^*(y,\tau)=u^*(x_0^+,t)\,,\\
&& u_*(x_0^-,t)={\rm ess}\!\!\!\!\!\!\!\!\liminf_{Q\ni(x,\tau)\to (x_0^-,t)}u_*(x,\tau)\,\leq G(x_0^-)+\ep+K(t-t_0)\,<\\
&&\qquad \qquad <\,G(x_0^+)-\ep+k(t-t_0)\,\leq\,{\rm ess}\!\!\!\!\!\!\!\!\liminf_{Q\ni(y,\tau)\to (x_0^+,t)}u^*(y,\tau) =u_*(x_0^+,t)
\end{eqnarray*}
(the equalities in these estimates follow from Lemma \ref{propis} applied to $(a,x_0)\times (0,T)$, respectively $(x_0,b)\times (0,T)$). 
Therefore it follows from Lemma \ref{propis} that
\begin{equation*}\label{abs*1}
u^*(x_0,t)={\rm ess}\!\!\!\!\!\!\limsup_{Q\ni(y,\tau)\to (x_0,t)}u^*(y,\tau)=u^*(x_0^+,t)>u^*(x_0^-,t)
\end{equation*}
\begin{equation*}\label{abs*2}
u_*(x_0,t)={\rm ess}\!\!\!\!\!\!\liminf_{Q\ni(x,\tau)\to (x_0,t)}u_*(x,\tau)=u_*(x_0^-,t)< u_*(x_0^+,t_0)
\end{equation*}
for all $t\in (t_0,\underline t_\ep)$. 
Since $\ep$ is arbitrary and $\underline t_\ep\to \underline t$ as $\ep\to0$, the conclusion follows. 
\end{proof}

The concept of {\it barrier effect of a discontinuity}, discussed in the Introduction, is made precise by the following lemma.
 
\begin{lem}\label{prerem} 
Let $-\infty<a<c<b<\infty$, and let $Q=(a,b)\times(0,T)$, $Q^-=(a,c)\times (0,T)$ and $Q^+=(c,b)\times (0,T)$. 
Let $u$ be a viscosity solution of problem $(N_{a,\pm}^{b,\pm})$.

\smallskip

\noindent $(i)$ If $u^*(c^+,\cdot)> u^*(c^-,\cdot)$ in  $(0,T)$, then $u\lefthalfcup Q^+$ is a viscosity solution of $(N_{c,+}^{b,\pm})$.
\smallskip 

\noindent $(ii)$ If $u_*(c^+,\cdot)> u_*(c^-,\cdot)$ in  $(0,T)$, then $u\lefthalfcup Q^-$ is a viscosity solution of $(N_{a,\pm}^{c,+})$.
\smallskip

\noindent $(iii)$ If $u^*(c^+,\cdot)< u^*(c^-,\cdot)$ in  $(0,T)$, then $u\lefthalfcup Q^-$ is a viscosity solution of  $(N_{a,\pm}^{c,-})$.
\smallskip

\noindent $(iv)$ If $u_*(c^+,\cdot)<u_*(c^-,\cdot)$ in  $(0,T)$, then $u\lefthalfcup Q^+$ is a viscosity solution of $(N_{c,-}^{b,\pm})$.
\smallskip

\noindent Similar statements hold for viscosity solutions of problems $(N_{a,\pm})$, $(N^{b,\pm})$ and equation \eqref{hje} in $\R$.
\end{lem}

\begin{remark}\label{rem.t12}
Let $0\leq t_1<t_2\leq T$ and $\Omega=(a,b)$, with $-\infty\leq a<b\leq \infty$. Then, for any $c\in (a,b)$, it can be easily checked that the conclusions of Lemma \ref{prerem} hold true for viscosity solutions of problems $(N_{a,\pm}^{b,\pm})$, $(N_{a,\pm})$ or $(N^{b,\pm})$ in $\Omega\times (t_1,t_2)$.
\end{remark}

\begin{proof} We only prove $(i)$, since the other proofs are similar. We set $\tilde u=u\lefthalfcup Q^+$. 
Since $\tilde u^*=u^*$ in $Q^+$  (see Remark \ref{restri}) and $u$ is a viscosity solution of $(N_{a,\pm}^{b,\pm})$ in $Q$,  
$\tilde u$  is a viscosity solution of \eqref{hje} in $Q^+$. By Definition \ref{defsubsupervi}$-(i)$,
it remains to prove that if $\varphi\in C^1([c,b]\times (0,T])$ and $(c,t_0)$ is a local maximum point of $\tilde u^*-\varphi$ in $[c,b]\times (0,T]$, then
 \begin{equation}\label{dimoc}
\varphi_t(c,t_0)+H(\varphi_x(c^+,t_0))\le 0 \qquad (t_0\in(0,T])\,.
\end{equation} 

 To prove \eqref{dimoc}, let $t_0\in (0,T)$ (if $t_0=T$ we argue as in \cite[Section 10.2]{Ev}) and observe first that, by assumption, there holds 
\begin{equation}\label{dimo1}
 u^*(c^-,t)<u^*(c,t)=u^*(c^+,t)=\tilde u^*(c,t)\quad \mbox{for all $t\in (0,T)$.}
\end{equation}
If $(c,t_0)$ is a strict local maximum point of $\tilde u^*-\varphi$ in $[c,b]\times (0,T]$, then 
\begin{equation}\label{con mez}
\tilde u^*(c,t_0)-\varphi(c,t_0)>\tilde u^*(y,\t)-\varphi(y,\t)\quad\text{for any $(y,\t)\in\overline{B_r^+(c,t_0)}$}
\end{equation}
for some $r>0$, where $B^+_r(c,t_0)=\{(y,\tau)\in B_r(c,t_0)\,|\,\,y>c\}$. 
Here $r$ is chosen such that $0<t_0-r<t_0+r<T$. In view of \eqref{dimo1} this implies that 
\begin{equation}\label{ugu}
u^*(y,\tau)=\tilde u^*(y,\tau)\ \  \mbox{for all $(y,\tau)\in \overline{B_r^+(c,t_0)}$}.
\end{equation} 
From \eqref{con mez} and \eqref{ugu} we obtain that
\begin{equation}\label{dimo2}
u^*(c,t_0)-\varphi(c,t_0)>u^*(y,\t)-\varphi(y,\t)\quad\text{for all $(y,\t)\in\overline{B_r^+(c,t_0)}\,.$}
\end{equation}
On the other hand, by \eqref{dimo1} and the upper semicontinuity of $u^*$ we also have that
$$
\limsup_{(y,\tau)\to (c,t), y<c} u^*(y,\tau)\le u^*(c^-,t)<u^*(c,t)\,,
$$
thus for some $r>0$ there holds
$$
u^*(c,t)>u^*(y,\t)\quad\text{for any $(y,\t)\in B_r^-(c,t_0)=\{(y,\tau)\in B_r(c,t_0)\,|\,\,y<c\}\,.$}
$$
Hence we can extend the definition of $\varphi$ to $[a,b]\times (0,T]$ so that $\varphi_x(c_-,t_0)=\varphi_x(c_+,t_0)$, and
\begin{equation}\label{dimo3}
u^*(c,t_0)-\varphi(c,t_0)>u^*(y,\t)-\varphi(y,\t)\quad\text{for any $(y,\t)\in B_r^-(c,t_0)\,.$}
\end{equation}
 By \eqref{dimo2}-\eqref{dimo3} $(c,t_0)$  is a local maximum point of $u^*-\varphi$ in $Q$, thus by \eqref{ad0-} we obtain \eqref{dimoc}.
\end{proof}  

To complete the proof of Proposition \ref{prolip}$-(ii)$, we show that $u^*-u_*$ is nonincreasing with respect to $t$ for viscosity solutions $u$ of 
equation \eqref{hje} and prove \eqref{desa1}.
\medskip

\noindent {\it Proposition \ref{prolip}\,$(ii)$: proof of \eqref{desa1}.}
We prove \eqref{desa+} assuming that $u_0(c^+)>u_0(c^-)$ and 
\begin{equation}\label{preme}
\begin{cases}
u^*(c,t)=u^*(c^+,t)>u^*(c^-,t)\\
u_*(c,t)=u_*(c^-,t)<u_*(c^+,t)
\end{cases}
\ \mbox{for all}\ \,t\in[0,t_c).
\end{equation}
Observe that \eqref{desa1} follows at once (by subtraction) if we show that
\begin{equation}\label{presa1} 
\frac{u^*(c,t_1)\!-\!u^*(c,t_0)}{t_1-t_0}\!\le\!-\limsup\limits_{s\to \infty}H(s),
\quad \frac{u_*(c,t_1)\!-\!u_*(c,t_0)}{t_1-t_0}\!\ge \!-\liminf\limits_{s\to \infty}H(s)
\end{equation}
for $0\le t_0<t_1<t_c$.

We only prove the first inequality of \eqref{presa1}. Let $b_0=b$ if $b<\infty$ and let $b_0\in (c,\infty)$ if $b=\infty$. 
We set $P=(c,b_0)\times (t_0,t_c)$ and $u_1=u\lefthalfcup P$.
By Remark \ref{restri} and \eqref{preme} 
\begin{equation}\label{preme2}
\begin{cases}
(u_1)^*(x,t)=u^*(x,t)&\mbox{for}\ \,(x,t)\in P\\
(u_1)^*(c,t)=u^*(c,t)&\mbox{for}\ \,t_0<t<t_c\,.
\end{cases}
\end{equation}
By Lemma \ref{propis} applied to $(u_1)^*$, for all $\ep>0$ there exists  $\delta>0$ such that
\begin{equation}\label{preme3}
(u_1)^*(x,t_0) \le (u_1)^*(c,t_0) +\ep      \quad\text{for all }x\in[c,c+\delta)\,.
\end{equation}
We set, for all $p>0$ and $(x,t)\in \overline{P}$,
$$
v(x,t;p)\,=\,\min\{ (u_1)^*(c,t_0)\!+\!\ep\!+\!p(x\!-\!c)\!-\!(t\!-\!t_0)H(p)\, , \, \|u\|_{L^{\infty}(P)}\!+\!K(t\!-\!t_0)\}.
$$
Observe that there exist $p_\ep\ge 0$ and $\xi\in (c,b_0)$ such that for all $p>p_\ep$ 
\begin{equation}\label{choice p 1}
v(x,t_0;p)\ge (u_1)^*(x,t_0) \quad\text{if }c\le x\le b_0
\end{equation}
and 
\begin{equation}\label{choice p 2}
v(x,t;p)=\|u\|_{L^{\infty}(P)}+K(t-t_0) \quad \text{if }(x,t)\in \overline{P}, \ x>\xi.
\end{equation}
One easily checks that if $p>p_\ep$, then, by \eqref{choice p 2}, $v(x,t;p)$ is a viscosity supersolution of problem $(N_{c,+}^{b_0,+})$ in $P$. 

On the other hand, if $b=b_0\in \R$ it follows from Lemma \ref{prerem}$-(i)$ (see also Remark \ref{rem.t12}) that $u_1$ is a solution of $(N_{c,+}^{b_0,\pm})$ in $P$, and hence, 
by Definition \ref{defsubsupervi}$-(i)$, $u_1$ is a viscosity subsolution of problem $(N_{c,+}^{b_0,+})$ in $P$ (if $b=\infty$ we argue similarly: the restriction $\tilde u=u \lefthalfcup \tilde P$, with $\tilde P=(c,\infty)\times (t_0,t_c)$, is a viscosity solution of $(N_{c,+})$ in $\tilde P$ and $u_1$, which coincides with the restriction of $\tilde u$ to $P$, is a viscosity subsolution of $(N_{c,+}^{b_0,+})$ in $P$).
By \eqref{choice p 1} and Theorem \ref{comp}, this implies that $(u_1)^*\leq v$ in $\overline{P}$ for all $p>p_{\ep}$. In particular 
$$
\begin{aligned}
u^*(c,t_1)&=(u_1)^*(c,t_1)\le (u_1)^*(c,t_0)+\ep -(t_1-t_0)H(p)\,\leq\\
&\leq \,u^*(c,t_0)+\ep -(t_1-t_0)H(p)\quad\text{for $t_0<t_1< t_c$}
\end{aligned}
$$
(here we have used the second equality in \eqref{preme2}).
Choosing $p=p_n$ such that $p_n\to +\infty$ and $H(p_n)\to \limsup\limits_{s\to \infty}H(s)$ as $n\to\infty$, we obtain that
$$
u^*(c,t_1)\le   u^*(c,t_0)+2\ep -(t_1-t_0)\limsup\limits_{s\to \infty}H(s)\,,
$$
and the first inequality in \eqref{presa1} follows from the arbitrariness of $\ep$. 
\hfill$\square$


\section{Proof of existence: Perron's method revisited}\label{perron}
\setcounter{equation}{0}

\noindent {\it Proof of Theorem \ref{exi}.} $(i)$ We first prove the existence of a viscosity solution of the Cauchy problem $(CP)$. We adapt Perron's method used by Ishii (\cite{I2}) to our definition of semi-continuous envelopes based on essential limits.

Let $u_0\in L^\infty_{\text{loc}}(\R)$, and let $K,\,k$ be defined in \eqref{Def K,k}. Set
$$
\underline u(x,t)= u_0(x)+kt, \quad \overline u(x,t)= u_0(x)+Kt \quad\text{for a.e.~}x\in \R,\ t\in [0,T].
$$
Then $\underline u$ is a viscosity subsolution of $(CP)$, $\overline u$ a viscosity supersolution, and
\begin{equation}\label{eq.f1}
\underline u^*(x,0)=\overline u^*(x,0)=(u_0)^*(x)\,, \quad \underline u_*(x,0)=\overline u_*(x,0)=(u_0)_*(x) \quad\text{for }x\in \R.
\end{equation}

Let $\mathcal S$ be the set of viscosity subsolutions $v$ of $(CP)$ such that  $v^*\le \overline u^*$ and $v_*\le \overline u_*$ in $\R\times [0,T]$ (observe that these two inequalities are equivalent: if for example $ v^*\le \overline u^*$, then $(v^*)_*\le (\overline u^*)_*$ and, by Lemma \ref{propis}, $ v_*\le \overline u_*$). We set 
\begin{equation}\label{deffu}
u(x,t)\,:=\, \sup \{v^*(x,t)\,|\, v\in \mathcal S\}\quad\text{for }x\in \R,\ t\in [0,T]\,.
\end{equation}
Since $\underline u\in \mathcal S$, we have $u\ge \underline u^*$  in $\R\times [0,T]$. In particular, it follows that $u^*\ge \underline u^*$ in $\R\times [0,T]$, hence  $u^*(x,0)\ge (u_0)^*(x)$ for $x\in \R$ (see \eqref{eq.f1}).  
On the other hand $u\le \overline u^*$, whence $u^*\le \overline u^*$ in $ \R\times [0,T]$ and, by \eqref{eq.f1}, $u^*(x,0)\le (u_0)^*(x)$  for $x\in \R$. So we have $u^*(x,0)= (u_0)^*(x)$. Analogously, since $u^*\ge \underline u^*$ and $u^*\le \overline u^*$ in $\R\times [0,T]$, from Lemma \ref{propis} we get $u_*=(u^*)_*\geq (\underline u^*)_*=\underline u_*$ and $u_*=(u^*)_*\leq (\overline u^*)_*=\overline u_*$ in $\R\times [0,T]$, thus (see \eqref{eq.f1}) 
$$
(u_0)_*(x)= \underline u_*(x,0)\leq u_*(x,0)\leq \overline u_*(x,0)=(u_0)_*(x)\quad\text{ for all $x\in\R$\,.}
$$ 
Then we have proved that $u_*(x,0)=(u_0)_*(x)$ for any $x\in\R$.

We claim that the function $u$ defined in \eqref{deffu} is a viscosity solution of $(CP)$. By the above remarks, it is enough to show that $u$ is a viscosity solution of $u_t+H(u_x)=0$ in $\R\times (0,T)$. We shall prove this in two steps.

\smallskip

\noindent{\em Step 1: $u$ is a viscosity subsolution of $u_t+H(u_x)=0$ in $\R\times (0,T)$.}

Let $v \in \mathcal S$, and fix $(x,t)\in \R\times (0,T)$. Then
\begin{equation}\label{u*>v^*}
\begin{aligned}
u^*(x,t)
&=\text{ess}\limsup_{(y,\tau)\to (x,t)} u(y,\tau)\,\ge\\
&\ge \text{ess}\limsup_{(y,\tau)\to (x,t)} v^*(y,\tau)
=\limsup_{(y,\tau)\to (x,t)} v^*(y,\tau)=v^*(x,t).
\end{aligned}
\end{equation}
Since there exists a sequence $\{v_n^*\}\subseteq \mathcal S$ such that $v_n^*(x,t)\to u(x,t)$, we have that  
$$
u^*(x,t)\ge  u(x,t).
$$
Since $(x,t)$ is arbitrary this statement holds for all $(x,t)\in \R\times (0,T)$.

Let $(x_0,t_0)\in \R\times (0,T)$ and $\varphi\in C^1(  \R\times (0,T))$ be such that $u^*-\varphi$ has a strict local maximum at $(x_0,t_0)$. Let 
$(y_n,\tau_n)\to (x_0,t_0)$ be such that 
$$
a_n\,:=\, u(y_n,\tau_n)-\varphi(y_n,\tau_n)\to u^*(x_0,t_0)-\varphi(x_0,t_0)
$$ 
(this is possible by the definition of $u^*$). By the definition of $u$, for all $n\in \N$ there exists $v_n\in \mathcal S$ such that 
$$
v_n^*(y_n,\tau_n)-\varphi(y_n,\tau_n)\ge a_n-\tfrac 1n.
$$
Let $r>0$ be so small that $u^*-\varphi<u^*(x_0,t_0)-\varphi(x_0,t_0)$ in $\overline B_r(x_0,t_0)\setminus \{(x_0,t_0)\} $, and let $(x_n,t_n)\in \overline B_r(x_0,t_0)$ be a point at which $v_n^*-\varphi$ has its maximum in $\overline B_r(x_0,t_0)$. Hence
$$
v_n^*(x_n,t_n)-\varphi(x_n,t_n)\ge v_n^*(y_n,\tau_n)-\varphi(y_n,\tau_n)\ge a_n-\tfrac 1n
$$
and, by \eqref{u*>v^*},
$$
v_n^*(x_n,t_n)-\varphi(x_n,t_n)\le u^*(x_n,t_n)-\varphi(x_n,t_n).
$$
This means that 
\begin{equation}\label{a_n}
u^*(x_n,t_n)-\varphi(x_n,t_n)\ge a_n-\tfrac 1n.
\end{equation}
Up to subsequences we have that $(x_n,t_n)\to (\overline x,\overline t)\in \overline B_r(x_0,t_0)$. Recalling that, by Lemma 2.1,
$$
u^*(\overline x,\overline t)=\limsup_{(x,t)\to (\overline x,\overline t)}u^*(x,t),
$$ 
for any $\varepsilon >0$ there exists $\delta>0$ such that 
$$
u^*(\overline x,\overline t)+\ep \ge \sup_{\overline B_\delta(\overline x,\overline t)}u^*.
$$ 
Since $(x_n,t_n)\to (\overline x,\overline t)$, for fixed $\delta$ there exists $N\in \N$ such that 
$$
(x_n,t_n)\in B_\delta(\overline x,\overline t)\quad\text{for }n>N.
$$
So we have found that for all $\ep>0$ there exists $N\in \N$ such that 
$$
u^*(\overline x,\overline t)+\ep\ge u^*(x_n,t_n)\quad\text{for }n>N.
$$
Combining this with \eqref{a_n} we find that
$$
u^*(\overline x,\overline t)-\varphi(x_n,t_n)\ge u^*(x_n,t_n)-\varphi(x_n,t_n)-\ep\ge  a_n-\tfrac 1n-\ep.
$$
Passing to the limit $n\to \infty$ and using the arbitrariness of $\ep$ we conclude that 
$$
u^*(\overline x,\overline t)-\varphi(\overline x,\overline t)\ge \lim_{n\to \infty} a_n=u^*(x_0,t_0)-\varphi(x_0,t_0).
$$
Since the local maximum of $u^*-\varphi$ at $(x_0,t_0)$ is strict, this means that $(\overline x,\overline t)=(x_0,t_0)$.
Finally, since $v_n^*-\varphi$ has a maximum at $(x_n,t_n)$ we have that 
$$
\varphi_t(x_n,t_n)+H(\varphi_x(x_n,t_n))\le 0.
$$
Passing to the limit $n\to\infty$ this implies that $\varphi_t(x_0,t_0)+H(\varphi_x(x_0,t_0))\le 0$, and we have completed Step 1.

\smallskip

\noindent{\em Step 2: $u$ is a viscosity solution of $u_t+H(u_x)=0$ in $\R\times (0,T)$.}

Arguing by contradiction we suppose the $u$ is not a viscosity solution of $u_t+H(u_x)=0$ in $\R\times (0,T)$.
In view of Step 1 this means that $u$ is not a viscosity supersolution of $u_t+H(u_x)=0$ in $\R\times (0,T)$.
Hence there exist $(x_0,t_0)\in \R\times (0,T)$ and $\varphi\in C^1(\R\times (0,T))$ such that $u_*-\varphi$ has a strict local minimum at $(x_0,t_0)$ and
$$
\varphi_t(x_0,t_0)+H(\varphi_x(x_0,t_0))< 0.
$$

Observe that by definition we have $u\leq \overline u^*$ in $\R\times [0,T]$, thus $u_*\leq (\overline u^*)_*=\overline u_*$ by Lemma \ref{prolip}. So let us assume that
\begin{equation}\label{eq.f2}
u_*(x_0,t_0)<\overline u_*(x_0,t_0)
\end{equation}
(otherwise also $\overline u_*-\varphi$ would have a minimum at $(x_0,t_0)$, and so  $\varphi_t+H(\varphi_x)\ge 0$ at $(x_0,t_0)$, a contradiction). Since $\varphi\in C^1$ and $\overline u_*$ is lower semicontinuous, there exist $r>0$ and $0<\delta_0\le \tfrac 12 [\overline u_*(x_0,t_0)-u_*(x_0,t_0)]$ such that 
$$
\begin{cases}
\varphi_t+H(\varphi_x)\le 0 \\
u_*(x_0,t_0)+\varphi-\varphi(x_0,t_0)+ \delta_0\le \overline u_*
\end{cases}
\quad\text{in }B_{2r}(x_0,t_0)\cap \big(\R\times [0,T]\big).
$$
Since the minimum of $u_*-\varphi$ at $(x_0,t_0)$ is strict and $u_*$ is lower semicontinuous, we may choose   $\delta\in (0,\delta_0]$ so small that  
$$
\text{$u_*(x_0,t_0)+\varphi-\varphi(x_0,t_0)+\delta<u_*$ in $B_{2r}(x_0,t_0)\setminus B_r(x_0,t_0)$.} 
$$
So 
$$
\begin{cases}
u_*(x_0,t_0)
+\varphi-\varphi(x_0,t_0)+ \delta_0\le \overline u^* &\text{in }B_{2r}(x_0,t_0)\cap \big(\R\times [0,T]\big)\\
u_*(x_0,t_0)
+\varphi-\varphi(x_0,t_0)+\delta<u^* &\text{in } B_{2r}(x_0,t_0)\setminus B_r(x_0,t_0).
\end{cases}
$$

Now we set
$$
w(x,t)\,:=\,
\begin{cases} 
\max\{u(x,t),u_*(x_0,t_0)+\varphi(x,t)-\varphi(x_0,t_0)+\delta\}&\text{if }(x,t)\in B_r(x_0,t_0)\\
u(x,t)&\text{otherwise}.
\end{cases}
$$
Since $\varphi_t+H(\varphi_x)\le 0$ in $B_{2r}(x_0,t_0)\cap \big(\R\times [0,T]\big)$, $w^*\le \overline u^*$ and $u\in \mathcal S$, the function $w$ belongs to $\mathcal S$.
This implies that for all $(x,t)\in B_r(x_0,t_0)$
$$
u(x,t)=\sup\{v^*(x,t)\,|\, v\in \mathcal S\}\ge w^*(x,t)\ge u_*(x_0,t_0)+\varphi(x,t)-\varphi(x_0,t_0)+\delta.
$$
Since the right-hand side of the above inequality is a smooth function, we also have that
$$
u_*(x,t)\ge  u_*(x_0,t_0)+\varphi(x,t)-\varphi(x_0,t_0)+\delta\quad \text{for }(x,t)\in B_r(x_0,t_0).
$$ 
Choosing $(x,t)=(x_0,t_0)$ we obtain $\delta\le 0$, a contradiction. This completes the proof of Theorem \ref{exi} in the case of problem $(CP)$.

\smallskip

\noindent $(ii)$  Let us now consider the initial-boundary value problem $(N)$. We shall only prove the result for problem $(N_{a,+}^{b,-})$ in $Q=(a,b)\times (0,T)$ with initial data $u_0\in L^{\infty}((a,b))$ $(-\infty<a<b<\infty);$ the proof is analogous in the other cases. 

Let $\tilde{u}_0\in L^{\infty}(\R)$ be defined by setting 
$$\tilde{u}_0(x)\,:=\,\left\{\begin{array}{ll}
u_0(x)&\mbox{for $a.e.$ $x\in (a,b)$}\,,\\
l_1&\mbox{for $x<a$}\,,\\
l_2&\mbox{for $x>b$}\,,
\end{array}\right.$$
where $l_1,l_2\,<\,0$ are chosen so that 
\begin{equation}\label{eq.li}
\left\{\begin{array}{ll}
-\,\displaystyle{\frac{\|u_0\|_{L^{\infty}((a,b))}+l_i}{K-k}\;>\,T}&\mbox{if $K>k$}\,,\\
\|u_0\|_{L^{\infty}((a,b))}+l_i\,<\,0&\mbox{if $K=k$}\,,
\end{array}\right.\qquad (i=1,2)
\end{equation}
and $K,k$ are the constants in \eqref{Def K,k}. Let $\tilde{u}$ be the viscosity solution of the Cauchy problem $(CP)$ with initial condition $\tilde{u}_0$, given by step $(i)$ above. By \eqref{cic} and \eqref{sih} for all $t\in (0,T)$ there holds 
\begin{eqnarray*}
&&(\tilde{u}_0)_*+kt=\tilde u_*(\cdot\,,0)+kt\leq 
\tilde{u}_*(\cdot\,,t)\leq \tilde{u}^*(\cdot\,,t)\leq\\
&\le& \tilde{u}^*(\cdot\,,0)+Kt=(\tilde{u}_0)^*+Kt\quad\mbox{in $\R$}
\end{eqnarray*}
(here $\tilde u_*\equiv (\tilde u)_*$ , $\tilde u^*\equiv (\tilde u)^*$  for notational simplicity). Moreover,
$$
\text{ $(\tilde{u}_0)^*(x)=l_1$ for $x<a$,  $(\tilde{u}_0)^*(x)=l_2$ for $x>b$, $(\tilde{u}_0)_*(x)\geq -\|u_0\|_{L^{\infty}((a,b))}$ for $x\in (a,b)$\,.}
$$ 
Then for all $t\in (0,T)$ we have  
$$
\tilde{u}^*(y,t)\leq l_1+Kt\quad\mbox{for $y<a$}\,,\quad \tilde{u}^*(y,t)\leq l_2+Kt\quad\mbox{for $y>b$}\,,
$$
$$
-\|u_0\|_{L^{\infty}((a,b))}+kt \leq \tilde{u}_*(x,t) \quad\mbox{for  $x\in (a,b)$}\,,
$$
whence by \eqref{eq.li}
\begin{equation*}
\tilde{u}^*(y,t)\leq l_1+Kt<-\|u_0\|_{L^{\infty}((a,b))}+kt\leq \tilde{u}_*(x,t)\quad\mbox{for all $y<a$ and $x\in (a,b)$}\,,
\end{equation*}
\begin{equation*}
\tilde{u}^*(y,t)\leq l_2+Kt<-\|u_0\|_{L^{\infty}((a,b))}+kt\leq \tilde{u}_*(x,t)\quad\mbox{for all $y>b$ and $x\in (a,b)$}\,.
\end{equation*}
Arguing as in the last part of the proof of Lemma \ref{ledi}, it follows from the above inequalities that  for all $t\in (0,T)$
$$
\tilde{u}^*(a^+,t)\geq \tilde{u}_*(a^+,t)>\tilde{u}^*(a^-,t), \quad \tilde{u}^*(b^-,t)\geq \tilde{u}_*(b^-,t)>\tilde{u}^*(b^+,t)\,.
$$ 
Hence by Lemma \ref{prerem}$-(i)$, $(iii)$ $u:=\tilde{u}\lefthalfcup Q$ is a viscosity solution of problem $(N_{a,+}^{b,-})$ with initial data $u_0$. This proves the result.
\hfill$\square$


\section{Proof of uniqueness}\label{unire}
\setcounter{equation}{0}

{In step $(ii)$ of the following proof we use  Lemma \ref{prerem}, which describes the barrier effect of spatial discontinuities, to handle possible discontinuities of viscosity solutions if $u_0$ is piecewise continuous. 

\smallskip

\noindent {\em Proof of Theorem \ref{unique}.}  The result is easily proven if $u_0\in C(\overline{\Omega})$: if $\Omega$ is bounded, then
$(u_0)^*=(u_0)_*=u_0$ in $\overline{\Omega}$ and, by \eqref{cic} and \eqref{cfrt}, $u^*=u_*= v^*=v_*$ in $\overline{Q}$;
hence  $u$ has a continuous representative $\tilde{u}$ in $\overline{Q}$ and, by Proposition \ref{prolip}, $\tilde{u}$ is  Lipschitz continuous with respect to $t$ in $\overline{Q}$ and satisfies \eqref{inelip}. If $\Omega$ is unbounded we argue similarly, using \eqref{cfrtbis}  instead of \eqref{cfrt}.

Let us consider the general case of initial data as in assumption $(H_2)$. For simplicity we suppose that $u_0$ has a single jump discontinuity at $x_1\in\Omega=(a,b)$, and that
\begin{equation}\label{uq1}
u_0(x_1^+)>u_0(x_1^-)\,.
\end{equation}
If $u_0(x_1^+)<u_0(x_1^-)$ or the number of jumps is finite, the proof is similar. 

Let $u$ and $v$ be two viscosity solutions of $(N)$ with initial datum $u_0$. By \eqref{uq1} and \eqref {wt 1} there exists $\underline{t}_1\in(0,T]$ such that
\begin{equation}\label{dimox}
t\in[0,\underline{t}_1)\;\Rightarrow\;
\begin{cases}
u^*(x_1,t)=u^*(x_1^+,t)> u^*(x_1^-,t), &u_*(x_1^+,t)> u_*(x_1^-,t)=u_*(x_1,t)\,,\\ 
v^*(x_1,t)=v^*(x_1^+,t)> v^*(x_1^-,t); &v_*(x_1^+,t)> v_*(x_1^-,t)=v_*(x_1,t)\,.
\end{cases}
\end{equation}
Therefore $\tau_1\,:=\,\sup\left\{\underline t_1\in(0,T)\,|\, \text{\eqref{dimox} holds}\right\}>0$. 

Consider the set $Q_{1,\tau_1}\,:=\,(a,x_1)\times (0,\tau_1)$ and the restrictions $u_1\,:=\,u\lefthalfcup Q_{1,\tau_1}$ and $v_1\,:=\,v\lefthalfcup Q_{j,\tau_1}$. By Lemma \ref{prerem} (see also Remark \ref{rem.t12}), $u_1$ and $v_1$ are viscosity solutions in $Q_{1,\tau_1}$ of a problem with lateral boundary condition $u_{1x}=v_{1x}=\infty$ at $x=x_1$, with continuous initial function $u_0\lefthalfcup (a,x_1)$.
As already observed, Theorem \ref{unique} holds for continuous initial data, hence there holds $u_1=v_1=\tilde u_1$ a.e.\ in $Q_{1,\tau_1}$, where $\tilde u_1\in C(\overline Q_{1,\tau_1})$.

Similarly, the restrictions of $u$ and $v$ to $Q_{2,\tau_1}\,:=\,(x_1,b)\times (0,\tau_1)$ coincide a.e.\ in $Q_{2,\tau_1}$ with a continuous solution 
$\tilde u_2\in C(\overline Q_{2,\tau_1})$. In particular we have found that  $u=v$ a.e.\ in $Q_{\tau_1}=(a,b)\times (0,\t_1)$, and \eqref{inelip} is
satisfied in $Q_{\tau_1}$.

It follows from \eqref{dimox} that for all $t\in (0,\tau_1)$
\begin{equation}\label{eq ide salti}
\begin{cases}
u^*(x_1,t)=u^*(x_1^+,t)=u_*(x_1^+,t)= \tilde{u}_2(x_1,t)\\
u_*(x_1,t)=u_*(x_1^-,t)=u^*(x_1^-,t)= \tilde{u}_1(x_1,t)
\end{cases}
\Rightarrow \ 
J(t):=\tilde{u}_2(x_1,t)-\tilde{u}_1(x_1,t)>0\,.
\end{equation}

If $\tau_1=T$ the proof is complete, otherwise we claim that $\tilde{u}_2(x_1,\tau_1)=\tilde{u}_1(x_1,\tau_1)$. 
Arguing by contradiction, it follows from the continuity of $\tilde{u}_j$ in $\overline{Q}_{j,\t_1}$ ($j=1,2$) 
that there exists $\eta>0$ such that $\tilde{u}_2(x_1,\tau)-\tilde{u}_1(x_1,\tau)\geq \eta$ for all $\tau\in(0,\tau_1)$ sufficiently close to $\tau_1$. Then by Lemma \ref{ledi} there exists $\delta>0$, independent of $\tau$ (see \eqref{cl21}), such that \eqref{dimox} holds for every $t\in (\tau,\tau+\delta)$, a contradiction for $\tau>\tau_1-\delta$ with the definition of $\tau_1$. 

Since $u^*=u_*=v^*=v_*=\tilde{u}_j$ in $Q_{j,\tau_1}$ and $\tilde{u}_j\in C(\overline{Q}_{j,\t_1})$, we have that \begin{equation}\label{uq3}
\sup_{x\in [a,b]}[u^*(x,t)-v_*(x,t)]_+\!=\!\sup_{x\in [a,b]}[v^*(x,t)-u_*(x,t)]_+ 
\!=\! J(t)\quad\text{for $t\in(0,\t_1)$}. 
\end{equation}
Let  $Q_{t,T}=(a,b)\times (t,T)$ $(t\in (0,\tau_1))$. For the sake of simplicity, we assume that $(a,b)$ is bounded (otherwise we argue similarly and consider suitable trapezoidal domains as in \eqref{trap bis} instead of $Q_{t,T}$). 
By Theorem \ref{comp} (applied in $Q_{t,T}$) and \eqref{uq3} 
\begin{equation*}
\sup_{\overline{Q}_{t,T}}\,[u^*-v_*]_+\leq \sup_{\overline{\Omega}}\,[u^*(x,t^+)-v_*(x,t^+)]_+\leq J(t)\,,
\end{equation*}
\begin{equation*}
\sup_{\overline{Q}_{t,T}}\,[v^*-u_*]_+\leq \sup_{\overline{\Omega}}\,[v^*(x,t^+)-u_*(x,t^+)]_+\leq J(t)\,.
\end{equation*}
Since $J(t)\to 0$ as $t\to \tau_1^-$, we conclude that $u^*=u_*=v^*=v_*$ in $\overline{Q}_{\t_1,T}$, whence $u=v$ $a.e.$ in $Q$.
Moreover, by the above considerations, both the restrictions $u_1=u\lefthalfcup Q_1$ and $u_2=u\lefthalfcup Q_2$ ($Q_1=(a,x_1)\times (0,T)$, $Q_2=(x_1,b)\times (0,T)$) of the viscosity solution $u$ of problem $(N)$ admit continuous representatives $\tilde u_i\in C(\overline{Q}_i)$ ($i=1,2$), and $\tilde u_1(x_1,t)\neq \tilde u_2(x_1,t)$ if and only if $t\in [0,\tau_1)$. This proves  claims $(a)$ and $(b)$ of Theorem \ref{unique}$(ii)$. Finally, \eqref{inelip} and claim $(d)$ follow from \eqref{sih} and \eqref{desa1}, since by \eqref{eq ide salti} we have $u^*(x_1,t)=\tilde u_2(x_1,t)$ and $u_*(x_1,t)=\tilde u_1(x_1,t)$ for all $t\in [0,\tau_1)$ (observe also that, if $\tau_1<T$, $u^*(x_1,t)=\tilde u_2(x_1,t)=u_*(x_1,t)=\tilde u_1(x_1,t)$ for $t\in [\tau_1,T]$).
\hfill$\square$
\medskip

Finally we show that the existence and uniqueness of piecewise continuous viscosity solutions implies a comparison principle for these solutions.

\smallskip

\noindent {\it Proof of Corollary \ref{thcs}}. 
Let $u_0\le v_0$ a.e.\ in $(a,b)$.
By the uniqueness of viscosity solutions (Theorem \ref{unique}), it is enough to show that the corresponding viscosity solutions $u$ and $v$ given by Perron's method
satisfy $u\le v$ a.e.\ in $(a,b)\times (0,T)$. But this follows trivially from the pointwise definition in the proof of Theorem \ref{exi} of $u$ and $v$ in terms of the 
corresponding sets $\mathcal S$, say $\mathcal S_u$ and $\mathcal S_v$, and the observation that $\mathcal S_u\subseteq \mathcal S_v$ since $u_0\le v_0$ a.e.\ in $(a,b)$.}
\hfill$\square$


\bigskip
\noindent
{\bf Acknowledgement.}
MB acknowledges the MIUR Excellence Department Project awarded to the Department of Mathematics, University of Rome Tor Vergata, CUP E83C18000100006,
and the project Beyond Borders of the University of Rome Tor Vergata, CUP E84I19002220005.



\end{document}